\RequirePackage{ifpdf}
\ifpdf 
\documentclass[pdftex]{sigma}
\else
\documentclass{sigma}
\fi

\newcommand{\mathsym}[1]{{}}

\begin{document}

\newcommand{\ac}[2]{a_{#1 #2}}
\newcommand{\al}{\alpha}
\newcommand{\Aut}{\operatorname{Aut}}
\newcommand{\bc}[2]{b_{#1 #2}}
\newcommand{\bt}{\mathbf{t}}
\newcommand{\bv}{\mathbf{v}}
\newcommand{\bw}{\mathbf{w}}
\newcommand{\eps}{\varepsilon}
\newcommand{\ga}{\gamma}
\newcommand{\Ga}{\Gamma}
\newcommand{\grad}{\operatorname{grad}}
\newcommand{\half}{{\tfrac{1}{2}}}
\newcommand{\hs}{\varphi}
\newcommand{\Id}{I\!d}
\newcommand{\M}{M^3}
\newcommand{\p}{\partial}
\newcommand{\pt}{\frac{\p}{\p t}}
\newcommand{\ptv}{\frac{\p}{\p v}}
\newcommand{\ptw}{\frac{\p}{\p w}}
\newcommand{\pp}[1]{\frac{\partial^2}{\partial {#1}^2}}
\newcommand{\R}{\mathbb{R}}
\newcommand{\Rf}{\R^4}
\newcommand{\ricLC}{\widehat{\operatorname{Ric}}}
\newcommand{\sign}{\operatorname{sign}}
\newcommand{\SO}{\mathrm{SO}}
\newcommand{\Span}{\mathrm{span}}
\newcommand{\trace}{\operatorname{trace}}
\newcommand{\Z}{\mathbb{Z}}

%
\newcounter{rom}
\renewcommand{\therom}{(\roman{rom})}
\newenvironment{romanlist}{\begin{list}{\therom}
                {\setlength{\leftmargin}{2em}\usecounter{rom}}}%
{\end{list}}
%

\allowdisplaybreaks

\renewcommand{\thefootnote}{$\star$}

\renewcommand{\PaperNumber}{097}

\FirstPageHeading

\ShortArticleName{Indef\/inite Af\/f\/ine Hyperspheres Admitting a Pointwise Symmetry. Part 2}

\ArticleName{Indef\/inite Af\/f\/ine Hyperspheres\\ Admitting a Pointwise Symmetry. Part 2\footnote{This paper is a
contribution to the Special Issue ``\'Elie Cartan and Dif\/ferential Geometry''. The
full collection is available at
\textit{}\href{http://www.emis.de/journals/SIGMA/Cartan.html}{http://www.emis.de/journals/SIGMA/Cartan.html}}}

\Author{Christine SCHARLACH}

\AuthorNameForHeading{C. Scharlach}

\Address{Technische Universit\"at Berlin,
Fak. II, Inst. f. Mathematik, MA 8-3, 10623 Berlin, Germany}
\Email{\href{mailto:schar@math.tu-berlin.de}{schar@math.tu-berlin.de}}
\URLaddress{\url{http://www.math.tu-berlin.de/~schar/}}

\ArticleDates{Received May 08, 2009, in f\/inal form October 06, 2009;  Published online October 19, 2009}

\Abstract{An af\/f\/ine hypersurface $M$ is said to admit a pointwise symmetry, if
there exists a subgroup $G$ of $\Aut(T_p M)$ for all $p\in M$, which
preserves (pointwise) the af\/f\/ine metric $h$, the dif\/ference tensor $K$
and the af\/f\/ine shape operator $S$. Here, we consider 3-dimensional
indef\/inite af\/f\/ine hyperspheres, i.e.\ $S= H\Id$ (and thus $S$ is
trivially preserved). In Part 1 we found the possible symmetry groups
$G$ and gave for each $G$ a canonical form of $K$. We started a
classif\/ication by showing that hyperspheres admitting a pointwise
$\Z_2\times \Z_2$ resp.\ $\R$-symmetry are well-known, they have
constant sectional curvature and Pick invariant $J<0$ resp.\
$J=0$. Here, we continue with af\/f\/ine hyperspheres admitting a
pointwise $\Z_3$- or $SO(2)$-symmetry. They turn out to be warped
products of af\/f\/ine spheres ($\Z_3$) or quadrics ($SO(2)$) with a
curve.}


\Classification{53A15; 53B30}

\renewcommand{\thefootnote}{\arabic{footnote}}
\setcounter{footnote}{0}

\newtheorem*{DecProb}{(De)composition Problem}

\section{Introduction}\label{sec:intro}

Let $M^n$ be a connected, oriented manifold. Consider an immersed
hypersurface with relative normalization, i.e., an immersion
$\hs\colon M^n \rightarrow \R^{n+1}$ together with a transverse vector
f\/ield $\xi$ such that $D \xi$ has its image in
$\hs_*T_xM$. Equi-af\/f\/ine geometry studies the properties of such
immersions under equi-af\/f\/ine transformations, i.e.\  volume-preserving
linear transformations ($SL(n+1,\R)$) and translations.

In the theory of nondegenerate equi-af\/f\/ine hypersurfaces there exists
a canonical choice of transverse vector f\/ield $\xi$ (unique up to
sign), called the af\/f\/ine (Blaschke) normal, which induces a connection
$\nabla$, a nondegenerate symmetric bilinear form $h$ and a 1-1 tensor
f\/ield $S$ by
\begin{gather}
	   D_X Y =\nabla_X Y +h(X,Y)\xi,\label{strGauss}\\
	   D_X \xi =-SX,\label{strWeingarten}
\end{gather}
for all $X,Y \in {\cal X}(M)$. The connection $\nabla$ is called the
induced af\/f\/ine connection, $h$ is called the af\/f\/ine metric or Blaschke
metric and $S$ is called the af\/f\/ine shape operator.  In general
$\nabla$ is not the Levi-Civita connection $\hat\nabla$ of $h$. The
dif\/ference tensor $K$ is def\/ined as
\begin{gather}\label{defK}
	  K(X,Y)=\nabla_X Y-\hat\nabla_X Y,
\end{gather}
for all $X,Y \in {\cal X}(M)$. Moreover the form $h(K(X,Y),Z)$ is a
symmetric cubic form with the property that for any f\/ixed $X\in {\cal
X}(M)$, $\trace K_X$ vanishes.  This last property is called the
apolarity condition. The dif\/ference tensor $K$, together with the
af\/f\/ine metric $h$ and the af\/f\/ine shape operator $S$ are the most
fundamental algebraic invariants for a nondegenerate af\/f\/ine
hypersurface (more details in Section~\ref{sec:basics}). We say that
$M^n$ is indef\/inite, def\/inite, etc.\ if the af\/f\/ine metric $h$ is
indef\/inite, def\/inite, etc.\ (Because the af\/f\/ine metric is a multiple of
the Euclidean second fundamental form, a positive def\/inite
hypersurface is locally strongly convex.) For details of the basic
theory of nondegenerate af\/f\/ine hypersurfaces we refer to \cite{LSZ93}
and~\cite{NS94}.
	
Here we will restrict ourselves to the case of af\/f\/ine hyperspheres,
i.e.\ the shape operator will be a (constant) multiple of the
identity ($S= H \Id$). Geometrically this means that all af\/f\/ine
normals pass through a f\/ixed point or they are parallel. There are many af\/f\/ine hyperspheres, even in the two-dimensional case only partial classif\/ications are known. This is due to the fact that af\/f\/ine hyperspheres reduce to the study of the Monge-Amp\`ere equations. Our question
is the following: {\em What can we say about a three-dimensional
affine hypersphere which admits a~pointwise $G$-symmetry, i.e.\ there
exists a non-trivial subgroup $G$ of the isometry group such that for
every $p\in M$ and every $L\in G$:}
\[
K(LX_p,L Y_p)= L(K(X_p, Y_p))\qquad \forall\,  X_p, Y_p\in T_p M.
\] We
have motivated this question in Part~1~\cite{S07a} (see also
\cite{Bry01,Vr04,LS05}). A classif\/ication of
$3$-dimen\-sional positive def\/inite af\/f\/ine hyperspheres admitting
pointwise symmetries was obtained in \cite{Vr04}. We continue the
classif\/ication in the indef\/inite case. We can assume that the af\/f\/ine
metric has index two, i.e.\ the corresponding isometry group is the
(special) Lorentz group $\SO(1,2)$. In Part 1, it turns out that a
$\SO(1,2)$-stabilizer of a nontrivial traceless cubic form is
isomorphic to either $\SO (2)$, $\SO(1,1)$, $\R$, the group $S_3$ of
order~6, $\Z_2 \times \Z_2$, $\Z_3$, $\Z_2$ or it is trivial. We have shown that hyperspheres admitting a pointwise $\Z_2 \times \Z_2$- resp.\ $\R$-symmetry are well-known, they have constant sectional curvature and Pick invariant $J<0$ resp.\ $J=0$.

In the following we classify the indef\/inite af\/f\/ine hyperspheres which
admit a pointwise $\Z_3$-, $\SO(2)$- or $\SO(1,1)$-symmetry. They turn out to
be warped
products of af\/f\/ine spheres ($\Z_3$) or quadrics ($SO(2)$, $SO(1,1)$) with a
curve. As a result we get a new composition method. Since the methods for the proofs for $SO(1,1)$ are similar to those for $\Z_3$- or $\SO(2)$-symmetry (and as long) we will omit them here. Both methods and results are dif\/ferent in case of $S_3$-symmetry and will be published elsewhere. The paper is organized as follows:

We will state the basic formulas of (equi-)af\/f\/ine hypersurface-theory needed in the further classif\/ication in Section~\ref{sec:basics}. In Section~\ref{sec:type2,3,4}, we show that in case of $\SO(2)$-, $S_3$- or $\Z_3$-symmetry we can extend the canonical form of $K$ (cf.~\cite{S07a}) locally.  Thus we can obtain information about the coef\/f\/icients of $K$ and $\nabla$ from the basic equations of Gauss,
Codazzi and Ricci (cf.~Section~\ref{sec:Intcond}). In Section~\ref{sec:type2,4} we show, that in case of $\Z_3$- or $\SO(2)$-symmetry it follows that the hypersurface
admits a warped product structure $\R\times_{e^f}N^2$. Then we
classify such hyperspheres by showing how they can be constructed
starting from $2$-dimensional positive def\/inite af\/f\/ine spheres
resp. quadrics (cf.\ Theorems~\ref{thm:ClassC1}--\ref{thm:ExC3}). We end in Section~\ref{sec:type8} by stating the classif\/ication results in case of $SO(1,1)$-symmetry (cf.\ Theorems~\ref{thm:ClassC1t8}--\ref{thm:ExC3t8}).

The classif\/ication can be seen as a generalization of the well
known Calabi product of hyperbolic af\/f\/ine spheres \cite{Ca72,HLV08}
and of the constructions for af\/f\/ine spheres considered in~\cite{DV94}. The following natural question for a (de)composition theorem,
related to these constructions, gives another motivation for studying $3$-dimensional
hypersurfaces admitting a pointwise symmetry:

\begin{DecProb} Let $M^n$ be a nondegenerate
affine hypersurface in $\R^{n+1}$. Under what conditions do there
exist affine hyperpsheres $M_1^r$ in $\R^{r+1}$ and $M_2^s$ in
$\R^{s+1}$, with $r+s=n-1$, such that $M = I \times_{f_1} M_1
\times_{f_2} M_2$, where $I \subset \R$ and $f_1$ and $f_2$ depend
only on $I$ $($i.e.~$M$ admits a warped product structure$)$? How can the
original immersion be recovered starting from the immersion of the
affine spheres?
\end{DecProb}

Of course the f\/irst dimension in which the above problem can be
considered is three and our results provide an
answer in that case.

\section{Basics of af\/f\/ine hypersphere theory}
\label{sec:basics}

First we recall the def\/inition of the af\/f\/ine normal $\xi$
(cf.~\cite{NS94}). In equi-af\/f\/ine hypersurface theory on the ambient
space $\R^{n+1}$ a f\/ixed volume form $\det$ is given. A transverse
vector f\/ield $\xi$ induces a volume form $\theta$ on $M$ by
$\theta(X_1,\ldots,X_n)=\det(\hs_* X_1,\ldots,\hs_* X_n,\xi)$. Also
the af\/f\/ine metric $h$ def\/ines a volume form $\omega_h$ on $M$, namely
$\omega_h=|\det h|^{1/2}$. Now the af\/f\/ine normal $\xi$ is uniquely
determined (up to sign) by the conditions that $D \xi$ is everywhere
tangential (which is equivalent to $\nabla \theta =0$) and that
\begin{gather}\label{BlaschkeNormal}
\theta = \omega_h.
\end{gather}
Since we only consider 3-dimensional indef\/inite hyperspheres, i.e.
\begin{gather}\label{affHypSphere}
S= H \Id,\qquad H=\text{const},
\end{gather}
we can f\/ix the orientation of the af\/f\/ine normal $\xi$ such that the
af\/f\/ine metric has signature one. Then the sign of $H$ in the
def\/inition of an af\/f\/ine hypersphere is an invariant.

Next we state some of the fundamental equations, which a
nondegenerate hypersurface has to satisfy, see also \cite{NS94} or
\cite{LSZ93}. These equations relate $S$ and $K$ with amongst others
the curvature tensor $R$ of the induced connection $\nabla$ and the
curvature tensor $\hat R$ of the Levi-Civita connection~$\widehat\nabla$ of the af\/f\/ine metric $h$.  There are the
Gauss equation for $\nabla$, which states that:
\begin{gather*}
R(X,Y)Z =h(Y,Z)SX -h(X,Z) SY,
\end{gather*}
and the Codazzi equation
\begin{gather*}
(\nabla_X S) Y =(\nabla_Y S) X.
\end{gather*}
Also we have the total symmetry of the af\/f\/ine cubic form
\begin{gather}\label{defC}
	C(X,Y,Z)= (\nabla_X h) (Y,Z) = -2 h(K(X,Y),Z).
\end{gather}
The fundamental existence and uniqueness theorem, see~\cite{Dil89} or~\cite{DNV90}, states that given~$h$,~$\nabla$ and~$S$ such that the
dif\/ference tensor is symmetric and traceless with respect to $h$, on a
simply connected manifold $M$ an af\/f\/ine immersion of $M$ exists if and
only if the above Gauss equation and Codazzi equation are satisf\/ied.

	From the Gauss equation and Codazzi equation above the Codazzi
equation for $K$ and the Gauss equation for $\widehat\nabla$ follow:
\begin{gather*}
(\widehat{\nabla}_X K)(Y,Z)- (\widehat{\nabla}_Y K)(X,Z)=  \half( h(Y,Z)SX-
h(X,Z)SY  - h(SY,Z)X + h(SX, Z)Y),\!
\end{gather*}
and
\begin{gather*} \hat R (X,Y)Z =
\tfrac{1}{2} (h(Y,Z) SX - h(X,Z) SY   + h(SY,Z)X - h(SX,Z)Y ) - [K_X
    ,K_Y] Z
\end{gather*}
If we def\/ine the Ricci tensor of the Levi-Civita connection $\widehat
\nabla$ by:
\begin{gather}\label{def:RicLC}
\ricLC(X,Y)=\trace\{ Z \mapsto \hat R(Z,X)Y\}.
\end{gather}
and the Pick invariant by:
\begin{gather}\label{def:Pick}
J= \frac{1}{n (n-1)} h(K,K),
\end{gather}
then from the Gauss equation we get for the scalar
curvature $\hat{\kappa}=\frac{1}{n (n-1)}(\sum_{i,j} h^{ij}
\ricLC_{ij})$:
\begin{gather}\label{TE}
\hat \kappa= H+J.
\end{gather}
For an af\/f\/ine hypersphere the Gauss and Codazzi equations have the form:
\begin{gather}\label{gaussInd}
 R(X,Y)Z=H(h(Y,Z)X -h(X,Z) Y),\\
\label{codazziS}  (\nabla_X H)Y=(\nabla_Y H)X, \qquad \text{i.e.}
\quad H= {\rm const},\\
\label{CodK}  (\widehat{\nabla}_X K)(Y,Z)=(\widehat{\nabla}_Y K)(X,Z),\\
\label{gaussLC}  \hat R (X,Y)Z = H(h(Y,Z)X -h(X,Z) Y)- [K_X ,K_Y] Z.
\end{gather}
Since $H$ is constant, we can rescale $\hs$ such that $H\in \{-1,0,1\}$.

\section[A local frame for pointwise $\SO(2)$-, $S_3$- or $\Z_3$-symmetry]{A local frame for pointwise $\boldsymbol{\SO(2)}$-, $\boldsymbol{S_3}$- or $\boldsymbol{\Z_3}$-symmetry}
\label{sec:type2,3,4}

Let $\M$ be a hypersphere admitting a $\SO(2)$-, $S_3$- or
$\Z_3$-symmetry. According to \cite[Theorem~4, 2.--4.]{S07a}  there exists
for every $p\in \M$ an ONB $\{ \bt, \bv, \bw\}$ of $T_p \M$ such that
\begin{alignat*}{4}
& K(\bt,\bt) = -2a_4 \bt, \qquad & & K(\bt,\bv) =a_4 \bv, \qquad && K(\bt,\bw) =
a_4 \bw, &\\
& K(\bv,\bv) = -a_4 \bt + a_6 \bv, \qquad && K(\bv,\bw) = -a_6 \bw, \qquad &&
K(\bw,\bw) = -a_4 \bt -a_6 \bv,&
\end{alignat*}
where $a_4 > 0$ and $a_6=0$ in case of $SO(2)$-symmetry, $a_4=0$ and
$a_6>0$ for $S_3$, and $a_4>0$ and $a_6>0$ for $\Z_3$.

We would like to extend the ONB locally. It is well known that
$\ricLC$ (cf.~\eqref{def:RicLC}) is a symmetric operator and we
compute
(some of the computations in this section are done with the
CAS
Mathematica\footnote{See Appendix or \url{http://www.math.tu-berlin.de/~schar/IndefSym_typ234.html}.}):
\begin{lemma} Let $p \in \M$ and $\{\bt,\bv,\bw\}$ the basis constructed
  earlier. Then
\begin{alignat*}{3}
  &\ricLC(\bt,\bt) =-2(H- 3a_4^2) , \qquad  &&\ricLC(\bt,\bv)=0, & \\
  &\ricLC(\bt,\bw)=0,\qquad  &&
\ricLC(\bv,\bv)=2(H-a_4^2+a_6^2) , & \\
  &\ricLC(\bv,\bw)=0,&&\ricLC(\bw,\bw)=2(H-a_4^2+a_6^2) .&
\end{alignat*}
\end{lemma}

\begin{proof}
The proof is a straight-forward computation using the Gauss
equation~\eqref{gaussLC}.  It follows e.g.\ that
\begin{gather*}
  \hat R(\bt,\bv)\bt = H \bv -K_{\bt}(a_4\bv)+K_{\bv}(-2 a_4 \bt) = H
  \bv -a_4^2 \bv -2 a_4^2\bv  = (H- 3 a_4^2)\bv,\\
  \hat
  R(\bt,\bw)\bt = H \bw -K_{\bt}(a_4 \bw) +K_{\bw}(-2 a_4 \bt) =H \bw
  - a_4^2 \bw -2 a_4^2 \bw = (H- 3 a_4^2)\bw,\\
  \hat
  R(\bt,\bv)\bw =-K_{\bt}(-a_6 \bw)+K_{\bv}(a_4\bw)=0.
\end{gather*}
From this it immediately follows that
\[
\ricLC(\bt,\bt) = -2(H-3a_4^2)
\]
and
\[
\ricLC(\bt,\bw)=0.
\]
The other equations follow by similar computations.
\end{proof}

We want to show that the basis we have constructed, at each point $p$,
can be extended dif\/ferentiably to a neighborhood of the point $p$ such
that, at every point, $K$ with respect to the frame $\{T,V,W\}$ has
the previously described form.
\begin{lemma}\label{lem:KfT234}
Let $\M$ be an affine hypersphere in $\mathbb R^4$ which admits a
  pointwise $\SO(2)$-, $S_3$- or $\Z_3$-symmetry. Let $p \in M$. Then
  there exists an orthonormal frame $\{T,V,W\}$ defined in a~neighborhood of the point $p$ such that $K$ is given by:
\begin{alignat*}{4}
& K(T,T)= -2a_4 T,\qquad && K(T,V) =a_4 V,\qquad && K(T,W) =
a_4 W,& \\
& K(V,V)= -a_4 T + a_6 V,\qquad  && K(V,W)= -a_6 W, \qquad &&
K(W,W) = -a_4 T -a_6 V, &
\end{alignat*}
where $a_4 > 0$ and $a_6=0$ in case of $\SO(2)$-symmetry, $a_4=0$ and
$a_6>0$ in case of $S_3$-symmetry, and $a_4>0$ and $a_6>0$ in case of
$\Z_3$-symmetry.
\end{lemma}

\begin{proof}First we want to show that at every point the vector
  $\bt$ is uniquely def\/ined (up to sign) and dif\/ferentiable. We
  introduce a symmetric operator $\hat A$ by:
\begin{gather*}
\ricLC(Y,Z)= h(\hat A Y,Z).
\end{gather*}
Clearly $\hat A$ is a dif\/ferentiable operator on $M$. Since $2(H-
3a_4^2) \neq 2(H-a_4^2+a_6^2)$, the operator has two distinct
eigenvalues. A standard result then implies that the
eigen distributions are dif\/ferentiable. We take $T$ a local unit
vector f\/ield spanning the 1-dimensional eigen distribution, and local
orthonormal vector f\/ields $\tilde{V}$ and $\tilde{W}$ spanning the
second eigen distribution. If $a_6=0$, we can take $V=\tilde{V}$ and
$W= \tilde{W}$.

As $T$ is (up to sign) uniquely determined, for $a_6\neq 0$ there
exist dif\/ferentiable functions~$a_4$,~$c_6$ and~$c_7$,
$c_6^2+c_7^2\neq 0$, such that
\begin{alignat*}{3}
& K(T,T)= -2a_4 T,\qquad & & K(\tilde{V},\tilde{V})= -a_4 T + c_6 \tilde{V}+
c_7 \tilde{W},& \\
& K(T,\tilde{V}) =a_4 \tilde{V},\qquad & &
K(\tilde{V},\tilde{W})= c_7 \tilde{V} -c_6 \tilde{W},& \\
& K(T,\tilde{W})= a_4 \tilde{W},\qquad && K(\tilde{W},\tilde{W})= -a_4 T -c_6
\tilde{V} -c_7 \tilde{W}.&
\end{alignat*}
As we have shown in \cite{S07a}, in
the proof of Theorem~2 (Case~2), we
can always rotate $\tilde{V}$ and $\tilde{W}$ such that we obtain the
desired frame.
\end{proof}

\begin{remark} It actually follows from the proof of the previous lemma
  that the vector f\/ield $T$ is (up to sign) invariantly def\/ined on
  $M$, and therefore the function $a_4$, too. Since the Pick invariant~\eqref{def:Pick} $J= \frac{1}{3}(-5 a_4^2 + 2 a_6^2)$, the function
  $a_6$ also is invariantly def\/ined on the af\/f\/ine hypersphere~$\M$.
\end{remark}

\section[Gauss and Codazzi for pointwise $\SO(2)$-, $S_3$- or $\Z_3$-symmetry]{Gauss and Codazzi for pointwise $\boldsymbol{\SO(2)}$-, $\boldsymbol{S_3}$- or $\boldsymbol{\Z_3}$-symmetry}
\label{sec:Intcond}

In this section we always will work with the local frame constructed
in the previous lemma. We denote the coef\/f\/icients of the Levi-Civita
connection with respect to this frame by:
\begin{alignat*}{4}
&   \widehat{\nabla}_T T = \ac12 V + \ac13 W,\qquad &&
   \widehat{\nabla}_T V = \ac12 T - \bc13 W,\qquad &&
   \widehat{\nabla}_T W = \ac13 T + \bc13 V,& \\
&   \widehat{\nabla}_V T = \ac22 V + \ac23 W, \qquad &&
   \widehat{\nabla}_V V = \ac22 T - \bc23 W, \qquad &&
   \widehat{\nabla}_V W = \ac23 T + \bc23 V,& \\
&   \widehat{\nabla}_W T = \ac32 V + \ac33 W, \qquad &&
   \widehat{\nabla}_W V = \ac32 T - \bc33 W, \qquad &&
   \widehat{\nabla}_W W = \ac33 T + \bc33 V.&
\end{alignat*}
We will evaluate f\/irst the Codazzi and then the Gauss equations
(\eqref{CodK} and \eqref{gaussLC}) to obtain more information.

\begin{lemma}\label{lem:CodK}
Let $\M$ be an affine hypersphere in $\mathbb R^4$ which admits a
pointwise $\SO(2)$-, $S_3$- or $\Z_3$-symmetry and $\{T,V,W\}$ the
corresponding ONB. If the symmetry group is
\begin{description}\itemsep=0pt
\item[$\mathbf{SO(2)}$,] then $0= \ac12 =\ac13 =\ac23 =\ac32$,
$\ac33=\ac22$ and \\ $T(a_4)=-4\ac22 a_4$, $0=V(a_4) = W(a_4)$,
\item[$\mathbf{S_3}$,] then $0= \ac12 =\ac13$, $ \ac23=-3 \bc13=
-\ac32$, $\ac33=\ac22$ and \\ $T(a_6)=-\ac22 a_6$, $V(a_6) = 3\bc33
a_6$, $W(a_6)=-3 \bc23 a_6$,
\item[$\mathbf{\Z_3}$ and $\mathbf{a_6^2\neq 4 a_4^2}$,] then $0=
\ac12 =\ac13 =\ac23 =\ac32$, $\ac33=\ac22$, $\bc13=0$,\\
$T(a_4)=-4\ac22 a_4$, $0=V(a_4) = W(a_4)$, and \\ $T(a_6)=-\ac22 a_6$,
$V(a_6) = 3\bc33 a_6$, $W(a_6)=-3 \bc23 a_6$,
\item[$\mathbf{\Z_3}$ and $\mathbf{a_6=2 a_4}$,] then $\ac12 =2 \ac22=
-2 \ac33=- \bc33$,\\ $\ac13 =- 2\ac23 = -2 \ac32= \bc23$, $\bc13=0$,
and \\ $T(a_4)=0$, $V(a_4)=-4 \ac22 a_4$, $W(a_4)= 4 \ac23 a_4$.
\end{description}
\end{lemma}

\begin{proof}
An evaluation of the Codazzi equations \eqref{CodK} with the help of
the CAS
Mathematica leads to the following equations (they relate
to eq1--eq6 and eq8--eq9 in the Mathematica notebook):
\begin{gather}
   V(a_4)=- 2\ac12 a_4, \qquad T(a_4)=-4 \ac22 a_4 + \ac12 a_6, \qquad
   0=4 \ac23 a_4 + \ac13 a_6,\label{CodKeq1} \\ W(a_4)=- 2\ac13
   a_4,\qquad 0=4 \ac32 a_4 + \ac13 a_6,\qquad T(a_4)=-4 \ac33 a_4 -
   \ac12 a_6, \label{CodKeq2}\\ T(a_6)-V(a_4)=3\ac12 a_4 -\ac22
   a_6,\qquad 0=\ac13 a_4 + (\ac23 + 3 \bc13)a_6, \label{CodKeq3}\\
   W(a_4)=(\ac23+ \ac32) a_6,\qquad  W(a_6)=
(-\ac23+3 \ac32)a_4 - \bc23 a_6,\nonumber \\ V(a_6)=(-\ac22+\ac33) a_4+ 3 \bc33
   a_6,
   \label{CodKeq4} \\
   T(a_6)=-\ac12 a_4 -\ac33 a_6,\qquad
   W(a_4)=- 3\ac13 a_4+ (-\ac32+3 \bc13)a_6,\label{CodKeq5}\\
   V(a_4)=(- \ac22+\ac33) a_6, \qquad W(a_6)= (3\ac23 - \ac32)a_4 -3
   \bc23 a_6,\label{CodKeq6}\\ 0=(\ac23-\ac32) a_4,\label{CodKeq8}\\
   W(a_4)=-\ac13 a_4 +(\ac32-3 \bc13) a_6.\label{CodKeq9}
 \end{gather}

 From the f\/irst equation of \eqref{CodKeq2} (we will use the notation
 \eqref{CodKeq2}.1) and \eqref{CodKeq4}.1 resp.\ \eqref{CodKeq1}.3 and~\eqref{CodKeq2}.2 we get:
 \begin{gather}\label{a13,a23+a32}
 0 =2 \ac13 a_4 +(\ac23+\ac32)a_6,\qquad
 0 =2 (\ac23+\ac32)a_4 + \ac13 a_6.
 \end{gather}
  From \eqref{CodKeq6}.1) and \eqref{CodKeq1}.1
  resp.~\eqref{CodKeq1}.2 and \eqref{CodKeq2}.3 we get:
 \begin{gather}\label{a12,a22+a33}
 0 =-2 \ac12 a_4 +2 (\ac22-\ac33)a_6,\qquad
 0 =2 (-\ac22+\ac33)a_4 + \ac12 a_6.
 \end{gather}
 We consider f\/irst the case, that $\mathbf{a_6^2\neq 4 a_4^2}$. Then
 we obtain from the foregoing equations that $\ac13=0$,
 $\ac32=-\ac23$, $\ac12=0$ and $\ac33=\ac22$. Furthermore it follows
 from \eqref{CodKeq1}.1 that $V(a_4)=0$, from \eqref{CodKeq1}.2 that
 $T(a_4)=-4 \ac22 a_4$ and from \eqref{CodKeq1}.3 that $\ac23
 a_4=0$. Equation \eqref{CodKeq2}.1 becomes $W(a_4)=0$, equation
 \eqref{CodKeq3}.2 $T(a_6)=-\ac22 a_6$ and \eqref{CodKeq3}.3 $(\ac23+3
 \bc13)a_6=0$. Finally equation \eqref{CodKeq4}.2 resp. 3 gives
 $W(a_6)=-3 \bc23 a_6$ and $V(a_6)= 3 \bc33 a_6$.

 In case of $SO(2)$-symmetry ($a_4>0$ and $a_6=0$) it follows that
 $\ac23=0$ and thus the statement of the theorem.

 In case of $S_3$-symmetry ($a_4=0$ and $a_6>0$) it follows that
 $\ac23=-3\bc13$ and thus the statement of the theorem.

 In case of $\Z_3$-symmetry ($a_4>0$ and $a_6>0$) it follows that
 $\ac23=0$ and $\bc13=0$ and thus the statement of the theorem.

 In case that $a_6= \pm 2 a_4$ ($\neq 0$), we can choose $V$, $W$ such
 that $\mathbf{a_6 = 2a_4}$. Now equations \eqref{CodKeq8},
 \eqref{CodKeq1}.3 and \eqref{CodKeq3}.3 lead to $\ac23=\ac32$,
 $\ac13=-2 \ac23$ and $\bc13=0$. A combination of \eqref{CodKeq1}.2
 and \eqref{CodKeq2}.3 gives $\ac12=(\ac22-\ac33)$, and then by
 equations \eqref{CodKeq3}.2, \eqref{CodKeq1}.1 and \eqref{CodKeq1}.2
 that $\ac33=-\ac22$. Thus $T(a_4)=0$ by \eqref{CodKeq1}.2,
 $V(a_4)=- 4 \ac22 a_4$ by \eqref{CodKeq1}.1 and $W(a_4)= 4 \ac22
 a_4$ by \eqref{CodKeq2}.1. Finally \eqref{CodKeq4}.2 and
 \eqref{CodKeq2}.1 resp. \eqref{CodKeq4}.3 and \eqref{CodKeq1}.1 imply
 that $\bc23=-\ac23$ resp. $\bc33=-\ac22$.
 \end{proof}

An evaluation of the Gauss equations \eqref{gaussLC} with the help of
the CAS Mathematica leads to the following:

\begin{lemma}\label{lem:GaussLC}
Let $\M$ be an affine hypersphere in $\mathbb R^4$ which admits a
pointwise $\SO(2)$-, $S_3$- or $\Z_3$-symmetry and $\{T,V,W\}$ the
corresponding ONB. Then
\begin{gather}
T(\ac22) = -\ac22^2 +\ac23^2 +H-3 a_4^2,\label{Gauss1.1}\\
T(\ac23) =-2\ac22\ac23,\label{Gauss1.2}\\ W(\ac22) + V(\ac23)
 =0,\label{Gauss1.3}\\ W(\ac23) - V(\ac22)  =0,\label{Gauss1.4}\\
V(\bc13) - T(\bc23) = \ac22\bc23 + (\ac23 +
\bc13)\bc33,\label{Gauss1.5}\\ T(\bc33) - W(\bc13) =
(\ac23+\bc13)\bc23 - \ac22\bc33 ,\label{Gauss1.6}\\ V(\bc33) -
W(\bc23) = -\ac22^2-\ac23^2 +2 \ac23\bc13 +\bc23^2 +\bc33^2 +H +a_4^2
+2 a_6^2. \label{Gauss1.7}
\end{gather}
If the symmetry group is $\Z_3$, then $a_6^2\neq 4 a_4^2$.
\end{lemma}

\begin{proof}
The equations relate to eq11--eq13 and eq16 in the Mathematica
notebook. If $a_6^2= 4 a_4^2 (\neq 0)$, then we obtain by equations
eq11.1 and eq12.3 resp. eq15.3 and eq12.3 that $2 V(\ac22)=-4 \ac22^2
-H+ 3 a_4^2$ resp. $2 W(\ac23)=4 \ac23^2+H-3 a_4^2$, thus $
V(\ac22)-W(\ac23)= -2 \ac22^2 -2 \ac23^2 -H +3 a_4^2$. This gives a
contradiction to eq13.3, namely $ V(\ac22)-W(\ac23)= -2 \ac22^2 -2
\ac23^2 -H -9 a_4^2$.
\end{proof}

\section[Pointwise $\Z_3$- or $\SO(2)$-symmetry]{Pointwise $\boldsymbol{\Z_3}$- or $\boldsymbol{\SO(2)}$-symmetry} \label{sec:type2,4}

The following methods only work in the case of $\Z_3$- or $\SO(2)$-symmetry,
therefore the case of $S_3$-symmetry will be considered elsewhere.
As the vector f\/ield $T$ is globally def\/ined, we can def\/ine the
distributions $L_1=\Span\{T\}$ and $L_2=\Span\{V,W\}$. In the
following we will investigate these distributions. For the terminology
we refer to \cite{Noe96}.
\begin{lemma}\label{L1}
The distribution $L_1$ is autoparallel $($totally geodesic$)$ with respect
to $\widehat\nabla$.
\end{lemma}
\begin{proof} From $\widehat{\nabla}_{T} T = \ac12 V + \ac13 W=0$
(cf.~Lemmas~\ref{lem:CodK} and~\ref{lem:GaussLC}) the claim follows immediately.
\end{proof}
\begin{lemma}\label{L2}
  The distribution $L_2$ is spherical with mean curvature normal
  $U_2=\ac22 T$.
\end{lemma}
\begin{proof} For $U_2=\ac22 T\in L_1=L_2^{\perp}$ we have
  $h(\widehat{\nabla}_{E_a} E_b, T)= h(E_a, E_b) h(U_2,T)$ for $E_a,
E_b\in \{V,W\}$, and $h(\widehat{\nabla}_{E_a} U_2, T)= h(E_a(\ac22) T
+ \ac22 \widehat{\nabla}_{E_a} T, T)=0$ (cf.\ Lemma~\ref{lem:CodK} and
\eqref{Gauss1.3}, \eqref{Gauss1.4}).
\end{proof}
\begin{remark} $\ac22$ is independent of the
  choice of ONB $\{V,W\}$. It therefore is a globally def\/ined
  function on $\M$.
\end{remark}
We introduce a coordinate function $t$ by $\pt :=T$. Using the
previous lemma, according to \cite{PR93}, we get:

\begin{lemma}\label{warped} $(\M,h)$ admits a warped product structure
  $\M=I \times_{e^f}N^2$ with $f: I \to \mathbb R$
  satisfying
\begin{gather}\label{deff}
\frac{\partial f}{\partial t}=\ac22.
\end{gather}
\end{lemma}
\begin{proof} Proposition~3 in \cite{PR93} gives the warped product structure with
warping function  $\lambda_2:I \to \mathbb R$. If we introduce $f=\ln
\lambda_2$, following the proof we see that $\ac22
T=U_2=-\grad(\ln\lambda_2)=-\grad f$.
\end{proof}

\begin{lemma}\label{lem:curvN2} The curvature of $N^2$ is
${}^NK(N^2)=e^{2f}(H+2 a_6^2 + a_4^2-\ac22^2)$.
\end{lemma}

\begin{proof} From Proposition~2 in \cite{PR93} we get the following relation
between the curvature tensor $\hat R$ of the warped product $\M$ and
the curvature tensor $\tilde R$ of the usual product of
pseudo-Riemannian manifolds ($X,Y,Z\in {\cal X}(M)$ resp. their
appropriate projections):
\begin{gather*} \hat R(X,Y)Z = \tilde{R} (X,Y)Z + h(Y,Z)(\widehat{\nabla}_X U_2 - h(X,U_2)U_2) - h(\widehat{\nabla}_X
U_2- h(X,U_2)U_2,Z)Y \\
\phantom{\hat R(X,Y)Z =}{} - h(X,Z)(\widehat{\nabla}_Y U_2 -
h(Y,U_2)U_2) + h(\widehat{\nabla}_Y U_2- h(Y,U_2)U_2,Z)X \\
\phantom{\hat R(X,Y)Z =}{} +h(U_2,U_2)(h(Y,Z)X-h(X,Z)Y).
\end{gather*}
Now $\tilde{R}(X,Y)Z={}^N\hat R(X,Y)Z$ for all $X,Y,Z\in TN^2$ and
otherwise zero (cf.\ \cite[page~89]{O'N83}, Corollary~58) and
$K(N^2)=K(V,W)= \frac{h(-\hat{R}(V,W)V,W)}{h(V,V)h(W,W)-h(V,W)^2}$
(cf.~\cite[page~77]{O'N83}, the curvature tensor has the opposite
sign). Since $h(X,Y)=e^{2f} {}^Nh(X,Y)$ for $X,Y\in TN^2$, it follows
that
\[
{}^NK(N^2)=e^{2f} h(-{}^N\hat{R}(V,W)V,W).
\] Finally we obtain by
the Gauss equation~\eqref{gaussLC} the last
ingredient for the computation: $\hat R(V,W)V$ $ =-(H+2 a_6^2 + a_4^2)
W$ (cf.~the Mathematica notebook).
\end{proof}

Summarized we have obtained the following structure equations
(cf.\ \eqref{strGauss}, \eqref{strWeingarten} and \eqref{defK}), where
$a_6=0$ in case of $\SO(2)$-symmetry resp. $\bc13=0$ in case of
$\Z_3$-symmetry:
\begin{gather}
D_T T = -2a_4 T- \xi, \label{D11}\\
D_T V = +a_4 V  - \bc13 W, \label{D12}\\
D_T W = +\bc13 V  + a_4 W, \label{D13}\\
D_V T =+(\ac22 +a_4) V,  \label{D21}\\
D_W T =+(\ac22 + a_4)W,  \label{D31}\\
D_V V =+ a_6 V -\bc23 W  +(\ac22 - a_4)T +\xi, \label{D22}\\
D_V W = +\bc23 V - a_6 W, \label{D23}\\
D_W V =-(\bc33 +a_6) W,  \label{D32}\\
D_W W = +(\bc33- a_6) V  +(\ac22 - a_4) T +\xi, \label{D33}
\\
\label{Dxi}
D_{X} \xi= -H X.
\end{gather}
The Codazzi and Gauss equations (\eqref{CodK} and \eqref{gaussLC})
have the form (cf.\ Lemmas~\ref{lem:CodK} and \ref{lem:GaussLC}):
\begin{gather}
 T(a_4) =-4\ac22 a_4, \qquad 0=V(a_4) = W(a_4),\label{Da4}\\
 T(a_6) =-\ac22 a_6,\qquad V(a_6) = 3\bc33 a_6,\qquad W(a_6)=-3 \bc23
 a_6,\label{Da6}\\
 T(\ac22) = -\ac22^2 +H-3 a_4^2, \qquad
 V(\ac22)=0, \qquad W(\ac22) =0,\label{Da22}\\
 V(\bc13) - T(\bc23) = \ac22\bc23 + \bc13\bc33,\label{Db1}\\
 T(\bc33) - W(\bc13) =\bc13\bc23 - \ac22\bc33 ,\label{Db2}\\
 V(\bc33) - W(\bc23) =
 -\ac22^2+\bc23^2 +\bc33^2 +H +a_4^2 +2 a_6^2,\label{Db3}
\end{gather}
where $a_6=0$ in case of $\SO(2)$-symmetry resp. $\bc13=0$ in case of
$\Z_3$-symmetry.

Our f\/irst goal is to f\/ind out how $N^2$ is immersed in $\Rf$, i.e.\ to
f\/ind an immersion independent of $t$. A look at the structure
equations \eqref{D11}--\eqref{Dxi} suggests to start with a linear
combination of~$T$ and~$\xi$.

We will solve the problem in two steps. First we look for a vector
f\/ield $X$ with $D_T X=\alpha X$ for some function~$\alpha$: We def\/ine
$X:=A T +\xi$ for some function $A$ on $\M$. Then $D_T X=\alpha X$ if\/f
$\alpha=-A$ and $\pt A= -A^2 +2a_4 A+ H$, and $A:=\ac22- a_4$ solves
the latter dif\/ferential equation. Next we want to multiply $X$ with
some function $\beta$ such that $D_T (\beta X)=0$: We def\/ine a positive
function $\beta$ on $\R$ as the solution of the dif\/ferential equation:
\begin{gather}\label{dtbeta}
\tfrac{\partial}{\partial t} \beta = (\ac22- a_4)\beta
\end{gather}
with initial condition $\beta(t_0)>0$. Then $D_T(\beta X)=0$ and by
\eqref{D21}, \eqref{Dxi} and \eqref{D31} we get (since
$\beta$, $\ac22$ and $a_4$ only depend on $t$):
\begin{gather}
D_{T}(\beta((\ac22- a_4)T +\xi)) =0,\label{eq31}\\
D_{V}(\beta((\ac22- a_4)T +\xi)) =\beta(\ac22^2-a_4^2-H)V ,\label{eq32}\\
D_{W}(\beta((\ac22- a_4)T +\xi)) =\beta(\ac22^2-a_4^2-H)W.\label{eq33}
\end{gather}
To obtain an immersion we need that $\nu:=\ac22^2-a_4^2-H$ vanishes
nowhere, but we only get:
\begin{lemma}\label{nu}
  The function $\nu=\ac22^2-a_4^2-H$ is globally defined,
  $\pt(e^{2f} \nu)=0$ and $\nu$ vanishes identically or nowhere on $\R$.
\end{lemma}

\begin{proof} Since $0=\pt {}^NK(N^2) = \pt(e^{2f}(2a_6^2-\nu))$
  (Lemma~\ref{lem:curvN2}) and $\pt(e^{2f}2 a_6^2)=0$ (cf.~\eqref{Da6} and~\eqref{deff}), we get that $\pt(e^{2f} \nu)=0$. Thus $\pt\nu=-2 (\pt
  f)\nu= -2\ac22\nu$.
\end{proof}

\subsection[The first case: $\nu \neq 0$ on $\M$]{The f\/irst case: $\boldsymbol{\nu \neq 0}$ on $\boldsymbol{\M}$}
\label{sec:case1}

We may, by translating $f$, i.e.\ by replacing $N^2$ with a homothetic
copy of itself, assume that $e^{2f} \nu =\eps_1$, where $\eps_1 =\pm 1$.

\begin{lemma}\label{defphi}
$\varPhi:=\beta ((\ac22- a_4)T +\xi)\colon
M^3 \to \R^4$ induces a proper affine sphere structure, say~$\tilde{\phi}$, mapping $N^2$ into a 3-dimensional linear subspace
of $\R^4$. $\tilde{\phi}$ is part of a quadric iff $a_6 =0$.
\end{lemma}
\begin{proof}
By \eqref{eq32} and \eqref{eq33} we have $\varPhi_*(E_a)= \beta
\nu E_a$ for $E_a\in \{V,W\}$. A further dif\/ferentiation, using \eqref{D22}
($\beta$ and $\nu$ only depend on $t$), gives:
\begin{gather*}
D_{V} \varPhi_*(V)  = \beta \nu D_{V} V = \beta \nu ((\ac22- a_4)T +a_6 V - \bc23 W +\xi)\\
\phantom{D_{V} \varPhi_*(V)}{} =a_6\varPhi_*(V)-\bc23 \varPhi_*(W) +\nu \varPhi
 =a_6\varPhi_*(V)-\bc23 \varPhi_*(W) +\eps_1 e^{-2f}\varPhi.
\end{gather*}
Similarly, we obtain the other derivatives, using \eqref{D23}--\eqref{D33}, thus:
\begin{gather}
D_{V} \varPhi_*(V) =   a_6\varPhi_*(V)  -\bc23 \varPhi_*(W) +
e^{-2f}\eps_1 \varPhi, \label{Dphivv}\\
D_{V} \varPhi_*(W) =  \bc23
\varPhi_*(V)  -a_6 \varPhi_*(W),  \label{Dphivw}\\
D_{W} \varPhi_*(V) =
 -(\bc33+ a_6)\varPhi_*(W),  \label{Dphiwv}\\
  D_{W} \varPhi_*(W) =
(\bc33-a_6)\varPhi_*(V)  +e^{-2f}\eps_1 \varPhi, \label{Dphiww}\\ D_{E_a}
\varPhi  =  \beta e^{-2f}\eps_1 E_a.\label{Dphiea}
\end{gather}
The foliation at $f=f_0$ gives an immersion of $N^2$ to $M^3$, say
$\pi_{f_0}$. Therefore, we can def\/ine an immersion of $N^2$ to $\R^4$
by $\tilde{\phi}:=\varPhi\circ\pi_{f_0}$, whose structure equations
are exactly the equations above when $f=f_0$. Hence, we know that
$\tilde{\phi}$ maps $N^2$ into
$\Span\{\varPhi_*(V),\varPhi_*(W),\varPhi\}$, an af\/f\/ine hyperplane of
$\R^4$ and $\pt\varPhi=0$ implies $\varPhi(t,v,w)=\tilde{\phi}(v,w)$.

We can read of\/f the coef\/f\/icients of the dif\/ference tensor
$K^{\tilde{\phi}}$ of $\tilde{\phi}$ (cf.~\eqref{strGauss} and
\eqref{defK}): $K^{\tilde{\phi}}(\tilde{V},\tilde{V})=a_6 \tilde{V}$,
$K^{\tilde{\phi}}(\tilde{V},\tilde{W})=-a_6
\tilde{W}$,$K^{\tilde{\phi}}(\tilde{W},\tilde{W})=- a_6 \tilde{V}$, and
see that $\trace (K^{\tilde{\phi}})_X$ vanishes. The af\/f\/ine metric
introduced by this immersion corresponds with the metric on $N^2$.
Thus $\eps_1 \tilde{\phi}$ is the af\/f\/ine normal of $\tilde{\phi}$ and
$\tilde{\phi}$ is a proper af\/f\/ine sphere with mean curvature
$\eps_1$.  Finally the vanishing of the dif\/ference tensor
characterizes quadrics.
\end{proof}
Our next goal is to f\/ind another linear combination of $T$ and $\xi$,
this time only depending on $t$. (Then we can express $T$ in terms
of $\phi$ and some function of $t$.)
\begin{lemma}\label{defdelta}
  Define $\delta := H T +(\ac22+a_4) \xi$. Then there exist a constant
  vector $C \in \R^4$ and a~function $a(t)$ such that
\[
 \delta(t)= a(t) C.
 \]
\end{lemma}
\begin{proof} Using \eqref{D21} resp. \eqref{D31} and
  \eqref{Dxi} we obtain that $D_{V}\delta = 0=D_{W} \delta$. Hence
  $\delta$ depends only on the variable $t$. Moreover, we get by
  \eqref{D11}, \eqref{Da22}, \eqref{Da4} and \eqref{Dxi} that
\begin{gather*}
  \pt\delta =D_{T} (H T+(\ac22+a_4)\xi)\\
\hphantom{\pt\delta}{} =H(-2a_4 T-\xi) + (-\ac22^2
  +H -3 a_4^2-4 \ac22 a_4)\xi -(\ac22 + a_4) H T\\
\hphantom{\pt\delta}{} =-(3 a_4+\ac22)(H
  T +(\ac22+a_4)\xi) =-(3 a_4+\ac22) \delta.
\end{gather*}
This implies that there exists a constant vector $C$ in $\Rf$ and a
function $a(t)$ such that $\delta(t)=a(t)C$.
\end{proof}
Notice that for an improper af\/f\/ine hypersphere ($H=0$) $\xi$ is
constant and parallel to $C$. Combining $\tilde{\phi}$ and $\delta$ we
obtain for $T$ (cf.\ Lemmas~\ref{defphi} and \ref{defdelta}) that
\begin{gather}\label{T}
T(t,v,w)= -\frac{a}{\nu}C +\frac{1}{\beta\nu}(\ac22+a_4)\tilde{\phi}(v,w).
\end{gather}

In the following we will use for the partial derivatives the
abbreviation $\hs_x:= \frac{\partial}{\partial x}\hs $, $x=t,v,w$.
\begin{lemma}\label{partialF}
\begin{gather*}
 \hs_t = -\frac{a}{\nu}C +\pt\left(\frac{1}{\beta \nu}\right)\tilde{\phi},\qquad
 \hs_v = \frac{1}{\beta\nu} \tilde{\phi}_v,\qquad
 \hs_w = \frac{1}{\beta\nu} \tilde{\phi}_w.
\end{gather*}
\end{lemma}
\begin{proof} As by \eqref{dtbeta} and Lemma~\ref{nu} $\pt
  \frac{1}{\beta\nu}= \frac{1}{\beta\nu}(\ac22+a_4)$, we obtain the
  equation for $\hs_t =T$ by~\eqref{T}. The other equations follow
  from \eqref{eq32} and \eqref{eq33}.
\end{proof}

It follows by the uniqueness theorem of f\/irst order dif\/ferential
equations and applying a~translation that we can write
\[
\hs(t,v,w)= \tilde{a}(t) C +\frac{1}{\beta\nu}(t)
\tilde{\phi}(v,w)
\] for a suitable function $\tilde{a}$ depending only
on the variable $t$. Since $C$ is transversal to the image of
$\tilde{\phi}$ (cf.~Lemmas~\ref{defphi} and \ref{defdelta},
$\nu\not\equiv 0$), we obtain that after applying an equiaf\/f\/ine
transformation we can write: $\hs(t,v,w) =(\gamma_1(t), \gamma_2(t)
\phi(v,w))$, in which $\tilde{\phi}(v,w)=(0,\phi(v,w))$.  Thus we have
proven the following:

\begin{theorem}\label{thm:ClassC1} Let $\M$ be an indefinite affine hypersphere
of~$\Rf$ which admits a pointwise $\mathbb
 Z_3$- or
  $SO(2)$-symmetry. Let $\ac22^2-a_4^2 \neq H$ for some $p\in
  \M$. Then
 $\M$ is affine equivalent to
\[
\hs:\ I\times N^2\to \Rf: \ (t,v,w)\mapsto (\gamma_1(t), \gamma_2(t)
  \phi(v,w)),
\]
 where $\phi: N^2 \to \mathbb R^3$ is a $($positive
  definite$)$ elliptic or hyperbolic affine sphere and $\gamma:I\to
  \mathbb R^2$ is a curve.  Moreover, if $\M$
  admits a pointwise $SO(2)$-symmetry then $N^2$ is either an
  ellipsoid or a two-sheeted hyperboloid.
\end{theorem}

We want to investigate the conditions imposed on the curve $\ga$. For
this we compute the derivatives of $\hs$:
\begin{alignat}{4}
 & \hs_t=(\ga_1',\ga_2' \phi ),\qquad & & \hs_v =(0,\ga_2\phi_v ),\qquad && \hs_w
=(0,\ga_2\phi_w ),& \nonumber\\
& \hs_{tt}=(\ga_1'',\ga_2'' \phi),\qquad &&
\hs_{tv}=(0,\ga_2'\phi_{v}),\qquad &&  \hs_{tw}=(0,\ga_2'\phi_w), & \label{DF} \\
& \hs_{vv}=(0,\ga_2\phi_{vv}),\qquad && \hs_{vw}=(0,\ga_2\phi_{vw}),\qquad &&
\hs_{ww}=(0,\ga_2'\phi_{ww}). &\nonumber
\end{alignat}
Furthermore we have to distinguish if $\M$ is proper ($H=\pm 1$) or
improper ($H=0$).

First we consider the case that $\M$ is proper, i.e. $\xi=-H\hs$. An
easy computation shows that the condition that $\xi$ is a transversal
vector f\/ield, namely $ 0\neq \det(\hs_t ,\hs_v, \hs_w, \xi)=-\ga_2^2
(\ga_1\ga_2' - \ga_1' \ga_2) \det(\phi_v, \phi_w, \phi)$, is
equivalent to $\ga_2\neq 0$ and $\ga_1\ga_2' - \ga_1' \ga_2 \neq
0$. To check the condition that $\xi$ is the Blaschke normal
(cf.~\eqref{BlaschkeNormal}), we need to compute the Blaschke metric
$h$, using \eqref{strGauss}, \eqref{DF},
\eqref{Dphivv}--\eqref{Dphiww} and the notation $r,s\in\{v,w\}$ and
$g$ for the Blaschke metric of $\phi$:
\begin{gather*}
\hs_{tt} = \cdots \hs_t +\frac{\ga_1'\ga_2'' - \ga_1''
\ga_2'}{H(\ga_1\ga_2' - \ga_1' \ga_2)} \xi,\qquad  \hs_{tr} = \text{tang}, \\
\hs_{rs}= \text{tang}  -\frac{\ga_1' \ga_2}{H(\ga_1\ga_2' - \ga_1'
\ga_2)}\eps_1 g\left(\frac{\p}{\p r}, \frac{\p}{\p s}\right)\xi.
\end{gather*}
We obtain that
\[
\det h = h_{tt}
(h_{vv}h_{ww}-h_{vw}^2)=\frac{\ga_1'\ga_2'' - \ga_1'' \ga_2'}{H^3
(\ga_1\ga_2' - \ga_1' \ga_2)^3}(\ga_1')^2 \ga_2^2 \det g.
\]
 Thus
\[
\ga_2^4 (\ga_1\ga_2' - \ga_1'
\ga_2)^2 \det(\phi_v,\phi_w,\phi)^2=\left|\frac{\ga_1'\ga_2'' - \ga_1''
\ga_2'}{(\ga_1\ga_2' - \ga_1' \ga_2)^3} (\ga_1')^2 \ga_2^2 \det
g\right|
\] is equivalent to  \eqref{BlaschkeNormal}. Since $\phi$ is a def\/inite proper af\/f\/ine sphere with normal
$-\eps_1\phi$, we can again use \eqref{BlaschkeNormal} to obtain
\[\xi=-H\hs \Longleftrightarrow \ga_2^2|\ga_1\ga_2' - \ga_1' \ga_2|^5=
 |\ga_1'\ga_2'' - \ga_1'' \ga_2'| (\ga_1')^2 \neq 0.\] From the
computations above ($g$ is positive def\/inite) also it follows that $\hs$
is indef\/inite if\/f either
\begin{gather*}
H\sign(\ga_1\ga_2' - \ga_1' \ga_2) =
\sign(\ga_1'\ga_2'' - \ga_1'' \ga_2')= \sign(\ga_1'\ga_2
\eps_1)
\qquad \mbox{or}\\
-H\sign(\ga_1\ga_2' - \ga_1' \ga_2) =
\sign(\ga_1'\ga_2'' - \ga_1'' \ga_2')= \sign(\ga_1'\ga_2
\eps_1).
\end{gather*}

Next we consider the case that $\M$ is improper, i.e.~$\xi$ is
constant. By Lemma~\ref{defdelta} $\xi$ is parallel to $C$ and thus
transversal to $\phi$. Hence we can apply an af\/f\/ine transformation to
obtain $\xi=(1,0,0,0)$. An easy computation shows that the condition
that $\xi$ is a transversal vector f\/ield, namely $0\neq \det(\hs_t
,\hs_v, \hs_w, \xi)=-\ga_2^2 \ga_2' \det(\phi_v, \phi_w, \phi)$, is
equivalent to $\ga_2\neq 0$ and $\ga_2' \neq 0$. To check the
condition that $\xi$ is the Blaschke normal
(cf.~\eqref{BlaschkeNormal}) we need to compute the Blaschke metric
$h$, using \eqref{strGauss}, \eqref{DF},
\eqref{Dphivv}--\eqref{Dphiww} and the notation $r,s\in\{v,w\}$ and
$g$ for the Blaschke metric of $\phi$:
\begin{gather*}
\hs_{tt} = \cdots \hs_t -\frac{\ga_1'\ga_2'' - \ga_1'' \ga_2'}{\ga_2'}
\xi,\qquad \hs_{tr}= \text{tang}, \qquad \hs_{rs}= \text{tang}  +\frac{\ga_1'
\ga_2}{\ga_2'}\eps_1 g\left(\frac{\p}{\p r}, \frac{\p}{\p s}\right)\xi.
\end{gather*}
We obtain that
\[
\det h = h_{tt}
(h_{vv}h_{ww}-h_{vw}^2)=-\frac{\ga_1'\ga_2'' - \ga_1''
\ga_2'}{(\ga_2')^3}(\ga_1')^2 \ga_2^2 \det g.
\]
 Thus
\eqref{BlaschkeNormal} is equivalent to
\[
\ga_2^4 (\ga_2')^2
\det(\phi_v,\phi_w,\phi)^2=\left|\frac{\ga_1'\ga_2'' - \ga_1''
\ga_2'}{(\ga_2')^3}(\ga_1')^2 \ga_2^2 \det g\right|.
\] Since $\phi$ is a
def\/inite proper af\/f\/ine sphere with normal $-\eps_1\phi$, we can again
use \eqref{BlaschkeNormal} to obtain
\[\xi=(1,0,0,0) \Longleftrightarrow \ga_2^2|\ga_2'|^5=
|\ga_1'\ga_2'' - \ga_1'' \ga_2'| (\ga_1')^2 \neq 0 .\] From the
computations above also it follows that $\hs$ is indef\/inite if\/f either
\begin{gather*}
- \sign(\ga_2')= \sign(\ga_1'\ga_2'' - \ga_1'' \ga_2')= \sign(\ga_1'
\ga_2 \eps_1) \qquad \mbox{or}\\
 \sign(\ga_2')= \sign(\ga_1'\ga_2'' - \ga_1''
\ga_2')= \sign(\ga_1' \ga_2 \eps_1).
\end{gather*}

So we have seen under which conditions we can construct a 3-dimensional indef\/inite af\/f\/ine hypersphere out of an af\/f\/ine sphere:
\begin{theorem}\label{thm:ConstrC1}
  Let $\phi:N^2 \to \mathbb R^3$ be a positive definite elliptic or
  hyperbolic affine sphere $($with mean curvature $\eps_1=\pm 1)$, and
  let $\gamma=(\ga_1,\ga_2): I \to \mathbb R^2$ be a curve. Define $\hs:I\times N^2\to \Rf$ by $\hs(t,v,w)= (\gamma_1(t), \gamma_2(t)
  \phi(v,w))$.
  \begin{enumerate}\itemsep=0pt
  \item[$(i)$] If $\ga$ satisfies $\ga_2^2|\ga_1\ga_2' - \ga_1'
  \ga_2|^5= \sign(\ga_1' \ga_2 \eps_1)(\ga_1'\ga_2'' - \ga_1'' \ga_2')
  (\ga_1')^2\neq 0$, then $\hs$ defines \mbox{a~$3$-dimensional} indefinite
  proper affine hypersphere.
  \item[$(ii)$] If $\ga$ satisfies $\ga_2^2|\ga_2'|^5=
  \sign(\ga_1' \ga_2 \eps_1)(\ga_1'\ga_2'' - \ga_1'' \ga_2')
  (\ga_1')^2\neq 0$, then $\hs$ defines a $3$-dimensional indefinite
  improper affine hypersphere.
  \end{enumerate}
\end{theorem}

Now we are ready to check the symmetries.
\begin{theorem}\label{thm:ExC1}
  Let $\phi:N^2 \to \mathbb R^3$ be a positive definite elliptic or
  hyperbolic affine sphere $($with mean curvature $\eps_1=\pm 1)$, and
  let $\gamma=(\ga_1,\ga_2): I \to \mathbb R^2$ be a curve such that
  $\hs(t,v,w)=(\gamma_1(t), \gamma_2(t) \phi(v,w ))$ defines a $3$-dimensional indefinite
  affine hypersphere. Then $\hs(N^2\times I)$  admits a pointwise $\mathbb Z_3$- or
 $SO(2)$-symmetry.
\end{theorem}

\begin{proof}
We already have shown that $\hs$ def\/ines a 3-dimensional indef\/inite
proper resp. improper af\/f\/ine hypersphere. To prove the symmetry we
need to compute~$K$. By assumption, $\phi$ is an af\/f\/ine sphere with
Blaschke normal $\xi^\phi =-\eps_1 \phi$. For the structure equations
\eqref{strGauss} we use the notation $\phi_{rs} = {}^\phi \Ga_{rs}^{u}
\phi_u - g_{rs} \eps_1 \phi$, $r,s,u \in \{v,w\}$. Furthermore we
introduce the notation $\al = \ga_1 \ga_2' - \ga_1' \ga_2$. Note that
$\al'=\ga_1 \ga_2'' - \ga_1'' \ga_2$. If $\hs$ is proper, using \eqref{DF},
we get the structure equations \eqref{strGauss} for $\hs$:
\begin{gather*}
\hs_{tt} = \frac{\al'}{\al} \hs_t +\frac{\ga_1'\ga_2'' - \ga_1''
\ga_2'}{H\al} \xi,\qquad \hs_{tr} = \frac{\ga_2'}{\ga_2} \hs_r, \\ \hs_{rs} =
{}^\phi \Ga_{rs}^{u} \hs_u - g_{rs} \eps_1 \frac{\ga_1 \ga_2}{\al} \hs_t -
g_{rs} \eps_1 \frac{\ga_1' \ga_2}{H\al}\xi.
\end{gather*}
We compute $K$ using \eqref{defC} and obtain:
\begin{gather*}
(\nabla_{\hs_t} h)(\hs_r,\hs_s)  = \left(\left(\frac{\ga_1 \ga_2}{\al}\right)'
\frac{\al}{\ga_1 \ga_2} -2 \frac{\ga_2'}{\ga_2}\right) h(\hs_r,\hs_s),\\
(\nabla_{\hs_r} h)(\hs_t,\hs_t)  = 0,
\end{gather*}
implying that $K_{\hs_t}$ restricted to the space spanned by $\hs_v$ and
$\hs_w$ is a multiple of the identity. Taking $T$ in direction of $\hs_t$,
we see that $\hs_v$ and $\hs_w$ are orthogonal to $T$. Thus we can
construct an ONB $\{T,V,W\}$ with $V,W$ spanning $\Span\{\hs_v,\hs_w\}$
such that $a_1 = 2 a_4$, $a_2=a_3=a_5=0$. By the considerations in
\cite[Section~4]{S07a} we see that $\hs$ admits a pointwise $\mathbb Z_3$- or
 $SO(2)$-symmetry. If~$\hs$ is improper, the proof runs completely analogous.
\end{proof}

\subsection[The second case: $\nu \equiv 0$ and $H\neq 0$ on $\M$]{The second case: $\boldsymbol{\nu \equiv 0}$ and $\boldsymbol{H\neq 0}$ on $\boldsymbol{\M}$}
\label{sec:case2}

Next, we consider the case that $H =\ac22^2 - a_4^2$ and $H\neq 0$ on
$\M$. It follows that $\ac22\neq \pm a_4$ on~$\M$.

We already have seen that $\M$ admits a warped product structure. The
map $\varPhi$ we have constructed in Lemma~\ref{defphi} will not def\/ine
an immersion (cf.~\eqref{eq32} and \eqref{eq33}). Anyhow, for a f\/ixed
point $t_0$, we get from \eqref{D22}--\eqref{D33}, \eqref{eq32} and
\eqref{eq33}, using the notation $\tilde{\xi}=(\ac22-a_4) T + \xi$:
\begin{gather*}
D_V V =a_6 V - \bc23 W +\tilde{\xi},\qquad
D_V W =\bc23 V - a_6 W,\\
D_W V =-(\bc33 + a_6) W,\qquad
D_W W =(\bc33 - a_6) V +\tilde{\xi},\qquad
D_{E_a}\tilde{\xi}=0, \qquad E_a\in \{V,W\}.
\end{gather*}
Thus, if $v$ and $w$ are local coordinates which span the second
distribution $L_2$, then we can interpret $\hs(t_0,v,w)$ as a
positive def\/inite improper af\/f\/ine sphere in a $3$-dimensional
linear subspace.

Moreover, we see that this improper af\/f\/ine sphere is a paraboloid
provided that $a_6(t_0, v,w)$ vanishes identically. From the
dif\/ferential equations \eqref{Da6} determining $a_6$, we see that this
is the case exactly when $a_6$ vanishes identically, i.e.\  when $\M$
admits a pointwise $SO(2)$-symmetry.

After applying a translation and a change of coordinates, we may
assume that
\begin{gather*}
\hs(t_0,v,w)=(v,w,f(v,w),0),
\end{gather*}
with af\/f\/ine normal $\tilde{\xi}(t_0,v,w)=(0,0,1,0)$. To obtain $T$
at $t_0$, we consider \eqref{D21} and \eqref{D31} and get that
\begin{gather*}
D_{E_a}(T-(\ac22+a_4) \hs) = 0,\qquad E_a, E_b\in \{V,W\}.
\end{gather*}

Evaluating at $t=t_0$, this means that there exists a constant vector
$C$, transversal to $\Span\{V,W,\xi\}$, such that
$T(t_0,v,w)=(\ac22+a_4)(t_0) \hs(t_0,v,w) +C$. Since $\ac22+a_4\neq 0$
every\-where, we can write:
\begin{gather}\label{Tt0}
  T(t_0,v,w)=\alpha_1 (v,w, f(v,w),\alpha_2),
\end{gather}
where $\alpha_1, \alpha_2\neq 0$ and we applied an equiaf\/f\/ine
transformation so that $C=(0,0,0,\alpha_1\alpha_2)$. To obtain
information about $D_T T$ we have that $D_T T= -2 a_4 T -\xi$
(cf.~\eqref{D11}) and $\xi=\tilde{\xi} - (\ac22- a_4)T$ by the
def\/inition of $\tilde{\xi}$.  Also we know that
$\tilde{\xi}(t_0,v,w)=(0,0,1,0)$ and by \eqref{eq31}--\eqref{eq33}
that $D_{X}(\beta\tilde{\xi})=0$, $X\in {\cal X}(M)$. Taking suitable
initial conditions for the function $\beta$ ($\beta(t_0)=1$), we get
that $\beta\tilde{\xi}=(0,0,1,0)$ and f\/inally the following vector
valued dif\/ferential equation:
\begin{gather*}
D_T T= (\ac22 -3 a_4) T -\frac1{\beta} (0,0,1,0).
\end{gather*}
Solving this dif\/ferential equation, taking into account the initial
conditions \eqref{Tt0} at $t=t_0$, we get that there exist functions
$\delta_1$ and $\delta_2$ depending only on $t$ such that
\begin{gather*}
  T(t,u,v)= (\delta_1(t) v,\delta_1(t) w, \delta_1(t) (f(v,w)
  +\delta_2(t)), \alpha_2 \delta_1(t)),
\end{gather*}
where $\delta_1(t_0)=\alpha_1$, $\delta_2(t_0)=0$, $\delta_1'(t)
=(\ac22 -3 a_4) \delta_1(t)$ and $\delta_2'(t) =\delta_1^{-1}(t)
\beta^{-1}(t)$.  As $T(t,v,w) =\tfrac{\partial \hs}{\partial
  t}(t,v,w)$ and $\hs(t_0,v,w) =(v,w,f(v,w),0)$ it follows by
integration that
\[
\hs(t,v,w)= (\gamma_1(t) v, \gamma_1(t) w, \gamma_1(t) f(v,w)
+\gamma_2(t) , \alpha_2 (\gamma_1(t)-1)),
\]
where $\gamma_1'(t)
=\delta_1(t)$, $\gamma_1(t_0)=1$, $\gamma_2(t_0)=0$ and $\gamma_2'(t)
=\delta_1(t)\delta_2(t)$.  After applying an af\/f\/ine transformation
we have shown:

\begin{theorem} \label{thm:ClassC2} Let $\M$ be an indefinite proper affine
  hypersphere of $\,\Rf$ which admits a pointwise $\mathbb Z_3$- or
 $SO(2)$-symmetry. Let
  $H= \ac22^2 -a_4^2(\neq 0)$ on $M^3$. Then $\M$ is affine equivalent with
\[
\hs: \ I\times N^2\to \Rf: \ (t,v,w)\mapsto (\gamma_1(t) v, \gamma_1(t)
  w, \gamma_1(t) f(v,w) +\gamma_2(t),\gamma_1(t)),
  \] where $\psi: N^2
  \to \mathbb R^3:(v,w) \mapsto (v,w,f(v,w))$ is a positive definite
  improper affine sphere with affine normal $(0,0,1)$ and $\gamma:I\to
  \mathbb R^2$ is a curve. Moreover, if $\M$ admits a
  pointwise $SO(2)$-symmetry then $N^2$ is an elliptic paraboloid.
\end{theorem}

We want to investigate the conditions imposed on the curve $\ga$. For
this we compute the derivatives of $\hs$:
\begin{gather}
\hs_t  =(\gamma_1' v,\gamma_1' w,\gamma_1' f(v,w)
+\gamma_2',\gamma_1'),\nonumber\\
\hs_v  =(\gamma_1,0,\gamma_1 f_v,0),\qquad
 \hs_w =(0,\gamma_1,\gamma_1 f_w,0),\nonumber\\
\hs_{tt} =(\gamma_1'' v, \gamma_1''w ,\gamma_1'' f(v,w) +
\gamma_2'',\gamma_1''),  \label{DF2}
\\
\hs_{tv} =\tfrac{\gamma_1'}{\gamma_1}
\hs_v,\qquad  \hs_{tw}=\tfrac{\gamma_1'}{\gamma_1} \hs_w,\nonumber\\
\hs_{vv} =(0,0,f_{vv}\gamma_1,0),\qquad
\hs_{vw}=(0,0,\gamma_1 f_{vw},0),\qquad  \hs_{ww}=(0,0,\gamma_1 f_{ww},0)
\nonumber.
\end{gather}

$\M$ is a proper hypersphere, i.e.\ $\xi=-H\hs$. An easy computation
shows that the condition that $\xi$ is a transversal vector f\/ield,
namely $0\neq \det(\hs_t ,\hs_v, \hs_w, \xi)=-H \ga_1^2 (\ga_1\ga_2' -
\ga_1'\ga_2)$, is equivalent to $\ga_1\neq 0$ and $\ga_1\ga_2' -
\ga_1'\ga_2\neq 0$. Since $(0,0,1,0)= \frac{\ga_1}{\ga_1\ga_2' -
\ga_1'\ga_2} \hs_t - \frac{\ga_1'}{\ga_1\ga_2' - \ga_1'\ga_2} \hs$, we
have the following structure equations:
\begin{gather}
\hs_{tt} =\left(\frac{\ga_1''}{\ga_1'} + \frac{\ga_1'\ga_2'' -
\ga_1''\ga_2'}{\ga_1'}\frac{\ga_1}{\ga_1\ga_2' - \ga_1'\ga_2}\right) \hs_t +
\frac{\ga_1'\ga_2'' - \ga_1''\ga_2'}{\ga_1\ga_2' -
\ga_1'\ga_2}\frac1H \xi, \nonumber \\
\hs_{tr} =\frac{\gamma_1'}{\gamma_1} \hs_r,\qquad
\hs_{rs} =\frac{\ga_1^2}{\ga_1\ga_2' - \ga_1'\ga_2} f_{rs} \hs_t +
\frac{\ga_1\ga_1'}{\ga_1\ga_2' - \ga_1'\ga_2} f_{rs} \frac1H \xi
.\label{streqC2}
\end{gather}

We obtain:
\[
\det h = h_{tt} (h_{vv}h_{ww}-h_{vw}^2)
=\frac{\ga_1'\ga_2'' - \ga_1'' \ga_2'}{H^3(\ga_1\ga_2' -
\ga_1'\ga_2)^3}\ga_1^2(\ga_1')^2 (f_{vv} f_{ww}-f_{vw}^2).
\]  Since
$\psi$ is a positive def\/inite improper af\/f\/ine sphere with af\/f\/ine
normal $(0,0,1)$, we get by \eqref{BlaschkeNormal} that $f_{vv}
f_{ww}-f_{vw}^2 =1$. Now \eqref{BlaschkeNormal} (for $\xi$) is
equivalent to
\[
\ga_1^4 (\ga_1\ga_2' -\ga_1'\ga_2)^2 =
\left|\frac{\ga_1'\ga_2'' - \ga_1'' \ga_2'}{(\ga_1\ga_2'
-\ga_1'\ga_2)^3}\right|\ga_1^2(\ga_1')^2 .
\] It follows that
\[\xi=-H\hs \Longleftrightarrow \ga_1^2|\ga_1\ga_2' -\ga_1'\ga_2|^5=
|\ga_1'\ga_2'' - \ga_1'' \ga_2'| (\ga_1')^2 \neq 0 .\] From the
computations above also it follows that $\hs$ is indef\/inite if\/f either
\begin{gather*}
 \sign(\ga_1'\ga_2'' - \ga_1'' \ga_2')= \sign(H(\ga_1\ga_2'
-\ga_1'\ga_2))= -\sign(\ga_1 \ga_1') \qquad \mbox{or}\\
 \sign(\ga_1'\ga_2'' -
\ga_1'' \ga_2')= -\sign(H(\ga_1\ga_2' -\ga_1'\ga_2))= -\sign(\ga_1
\ga_1').
\end{gather*}

So we have seen under which conditions we can construct a 3-dimensional indef\/inite af\/f\/ine hypersphere out of an af\/f\/ine sphere:
\begin{theorem}\label{thm:ConstrC2}
  Let $\psi: N^2 \to \mathbb R^3:(v,w) \mapsto (v,w,f(v,w))$ be a
  positive definite improper affine sphere with affine normal
  $(0,0,1)$, and let $\gamma: I \to \mathbb R^2$ be a curve. Define $\hs:I\times N^2\to \Rf$ by $\hs(t,v,w)= (\gamma_1(t) v, \gamma_1(t) w, \gamma_1(t) f(v,w)
  +\gamma_2(t),\gamma_1(t))$.
  If $\ga=(\ga_1,\ga_2)$ satisfies
  $\ga_1^2|\ga_1\ga_2' -\ga_1'\ga_2|^5=
  -\sign(\ga_1\ga_1')(\ga_1'\ga_2'' - \ga_1'' \ga_2') (\ga_1')^2 \neq
  0$, then $\hs$  defines a $3$-dimensional indefinite proper affine
  hypersphere.
\end{theorem}

Now we are ready to check the symmetries.
\begin{theorem}\label{thm:ExC2}
  Let $\psi: N^2 \to \mathbb R^3:(v,w) \mapsto (v,w,f(v,w))$ be a
  positive definite improper affine sphere with affine normal
  $(0,0,1)$, and let $\gamma: I \to \mathbb R^2$ be a curve such that
  $\hs(t,v,w)=(\gamma_1(t) v, \gamma_1(t) w, \gamma_1(t) f(v,w)
  +\gamma_2(t),\gamma_1(t))$ defines a $3$-dimensional indefinite proper affine
  hypersphere. Then $\hs(N^2\times I)$ admits a pointwise $\mathbb
  Z_3$- or $SO(2)$-symmetry.
\end{theorem}
\begin{proof}
We already have shown that $\hs$ def\/ines a 3-dimensional indef\/inite
proper af\/f\/ine hypersphere with af\/f\/ine normal $\xi= -H\hs$. To prove the
symmetry we need to compute $K$. We get the induced connection and the
af\/f\/ine metric from the structure equations \eqref{streqC2}. We
compute~$K$ using \eqref{defC} and obtain:
\begin{gather*}
(\nabla_{\hs_t} h)(\hs_r,\hs_s) = \left(\pt\ln\left(\frac{\ga_1
\ga_1'}{\ga_1\ga_2' -\ga_1'\ga_2}\right)-2 \frac{\ga_1'}{\ga_1}\right)
h(\hs_r,\hs_s),\\ (\nabla_{\hs_r} h)(\hs_t,\hs_t) = 0,
\end{gather*}
implying that $K_{\hs_t}$ restricted to the space spanned by $\hs_v$ and
$\hs_w$ is a multiple of the identity. Taking $T$ in direction of $\hs_t$,
we see that $\hs_v$ and $\hs_w$ are orthogonal to $T$. Thus we can
construct an ONB $\{T,V,W\}$ with $V,W$ spanning $\Span\{\hs_v,\hs_w\}$
such that $a_1 = 2 a_4$, $a_2=a_3=a_5=0$. By the considerations in
\cite[Section~4]{S07a} we see that $\hs$ admits a pointwise
$\mathbb Z_3$- or $SO(2)$-symmetry.
\end{proof}

\subsection[The third case: $\nu \equiv 0$ and $H=0$ on $\M$]{The third case: $\boldsymbol{\nu \equiv 0}$ and $\boldsymbol{H=0}$ on $\boldsymbol{\M}$}
\label{sec:case3}

 The f\/inal cases now are that $\nu \equiv 0$ and $H=0$ on the
whole of $M^3$ and hence $\ac22=\pm a_4$.

First we consider the case that $\ac22= a_4 =:a>0$. Again we use that
$M^3$ admits a warped product structure and we f\/ix a parameter
$t_0$. At the point $t_0$, we have by \eqref{D22}--\eqref{Dxi}:
\begin{gather*}
D_V V =+ a_6 V -\bc23 W   +\xi, \\
D_V W = +\bc23 V - a_6 W,\\
D_W V =-(\bc33 +a_6) W, \\
D_W W = +(\bc33- a_6) V  +\xi,\\
D_{X} \xi= 0.
\end{gather*}

Thus, if $v$ and $w$ are local coordinates which span the second
distribution $L_2$, then we can interpret $\hs(t_0,v,w)$ as a
positive def\/inite improper af\/f\/ine sphere in a $3$-dimensional
linear subspace.

Moreover, we see that this improper af\/f\/ine sphere is a paraboloid
provided that $a_6(t_0, v,w)$ vanishes identically. From the
dif\/ferential equations \eqref{Da6} determining $a_6$, we see that this
is the case exactly when $a_6$ vanishes identically, i.e.\  when $\M$
admits a pointwise $SO(2)$-symmetry.

After applying an af\/f\/ine transformation and a change of coordinates,
we may assume that
\begin{gather}\label{initcond1}
\hs(t_0,v,w)=(v,w,f(v,w),0),
\end{gather}
with af\/f\/ine normal $\xi(t_0,v,w)=(0,0,1,0)$, actually
\[
\xi(t,v,w)=(0,0,1,0)
\] ($\xi$ is constant on $\M$ by
assumption). Furthermore we obtain by \eqref{D21} and \eqref{D31},
that $D_U T= 2 a U$ for all $U\in L_2$. We def\/ine $\delta:= T- 2a
\hs$, which is transversal to $\Span\{V,W, \xi\}$. Since $a$ is
independent of $v$ and $w$ (cf.~\eqref{Da4}), $D_U \delta =0$, and we
can assume that
\begin{gather}\label{initcond2}
T(t_0,v,w) - 2a(t_0) \hs(t_0,v,w)= (0,0,0,1).
\end{gather}
We can integrate \eqref{Da4} ($T(a)= -4 a^2$) and we take $a=
\frac1{4t}$, $t>0$. Thus \eqref{D11} becomes $D_T T=-\frac1{2t} T
-\xi$ and we obtain the following linear second order ordinary
dif\/ferential equation:
\begin{gather}\label{DTT}
\pp{t}\hs + \frac1{2t} \pt \hs = -\xi.
\end{gather}
The general solution is $\hs(t,v,w)= -\frac{t^2}3 \xi + 2\sqrt{t}
A(v,w)+ B(v,w)$. The initial conditions \eqref{initcond1} and
\eqref{initcond2} imply that $A(v,w)=\big(\frac{v}{2\sqrt{t_0}},
\frac{w}{2\sqrt{t_0}}, \frac{f(v,w)}{2\sqrt{t_0}}+\frac23 t_0^{3/2},
\sqrt{t_0}\big)$ and $B(v,w)=(0,0,-t_0^2,-2 t_0)$. Obviously we can
translate $B$ to zero. Furthermore we can translate the af\/f\/ine sphere
and apply an af\/f\/ine transformation to obtain
$A(v,w)=\frac{1}{2\sqrt{t_0}}(v, w, f(v,w), 1)$. After a change of
coordinates we get:
\begin{gather}\label{result3.1}
\hs(t,v,w)= (tv, tw, t f(v,w) - c t^4, t), \qquad c,t>0.
\end{gather}
Next we consider the case that $-\ac22= a_4 =:a>0$. Again we use that
$M^3$ admits a warped product structure and we f\/ix a parameter
$t_0$. A look at \eqref{D22}--\eqref{Dxi} suggests to def\/ine $\tilde{\xi}= -2a T+\xi$, then we get at the point $t_0$:
\begin{gather*}
D_V V =+ a_6 V -\bc23 W   +\tilde{\xi}, \\
D_V W = +\bc23 V - a_6 W,\\
D_W V =-(\bc33 +a_6) W, \\
D_W W = +(\bc33- a_6) V  +\tilde{\xi},\\
D_{U} \tilde{\xi}= 0.
\end{gather*}

Thus, if $v$ and $w$ are local coordinates which span the second
distribution $L_2$, then we can interpret $\hs(t_0,v,w)$ as a
positive def\/inite improper af\/f\/ine sphere in a $3$-dimensional
linear subspace.

Moreover, we see that this improper af\/f\/ine sphere is a paraboloid
provided that $a_6(t_0, v,w)$ vanishes identically. From the
dif\/ferential equations \eqref{Da6} determining $a_6$, we see that this
is the case exactly when $a_6$ vanishes identically, i.e.\  when $\M$
admits a pointwise $SO(2)$-symmetry.

After applying an af\/f\/ine transformation and a change of coordinates,
we may assume that
\begin{gather}\label{init}
\hs(t_0,v,w)=(v,w,f(v,w),0),
\end{gather}
with af\/f\/ine normal
\begin{gather}\label{initxi}
\tilde{\xi}(t_0,v,w)=(0,0,1,0).
\end{gather}
We have considered $\tilde{\xi}$ before. We can solve
\eqref{dtbeta} ($\pt \beta= -2a \beta$) explicitly by $\beta=c
\frac1{\sqrt{a}}$ (cf.~\eqref{Da4}) and get by
\eqref{eq31}--\eqref{eq33} that $D_X (\frac1{\sqrt{a}} \tilde{\xi})
=0$. Thus $\frac1{\sqrt{a}}(-2a T+\xi)=:C$ for a constant vector $C$,
i.e.  $T=-\frac1{2a}(\sqrt{a} C - \xi)$. Notice that by \eqref{Dxi}
$\xi$ is a constant vector, too. We can choose $a=\frac1{4|t|}$, $t<0$
(cf.~\eqref{Da4}), and we obtain the ordinary dif\/ferential equation:
\begin{gather}\label{DT}
\pt \hs= -\sqrt{|t|} C -2t \xi, \qquad t<0.
\end{gather}
 The solution (after a
translation) with respect to the initial condition \eqref{init} is
$\hs(t,v,w)=\frac23 |t|^{\frac32} C- t^2 \xi + (v,w,f(v,w),0)$. Notice
that $C$ is a multiple of $\tilde{\xi}$ and hence by \eqref{initxi} a
constant multiple of $(0,0,1,0)$. Furthermore $\xi$ is transversal to
the space spanned by
$\hs(t_0,v,w)$. So we get after an af\/f\/ine transformation and a change
of coordinates:
\begin{gather}\label{result3.2}
\hs(t,v,w)= (v,w, f(v,w) + c t^3, t^4 ),\qquad c,t>0.
\end{gather}
Combining both results \eqref{result3.1} and \eqref{result3.2} we have:

\begin{theorem} \label{thm:ClassC3} Let $\M$ be an indefinite improper affine
  hypersphere of $\Rf$ which admits a pointwise $\mathbb
  Z_3$- or $SO(2)$-symmetry. Let $\ac22^2 =a_4^2$ on $M^3$. Then $\M$ is affine
  equivalent with either
  \begin{gather*} \hs: \ I\times N^2\to \Rf: \ (t,v,w)\mapsto (t v, t w,t f(v,w) -c
  t^4,t),\qquad (\ac22 =a_4)\qquad \text{or}\\
   \hs: \ I\times N^2\to \Rf: \ (t,v,w)\mapsto (v, w, f(v,w) +c
  t^3,t^4),\qquad(-\ac22 =a_4),
\end{gather*}
 where $\psi: N^2 \to \mathbb R^3:(v,w) \mapsto (v,w,f(v,w))$ is a
  positive definite improper affine sphere with affine normal
  $(0,0,1)$ and $c,t\in \R^+$.
Moreover, if $\M$ admits a pointwise $SO(2)$-symmetry then~$N^2$ is an
  elliptic paraboloid.
\end{theorem}

The computations for the converse statement can be done completely
analogous to the previous cases, they even are simpler (the curve is
given parametrized).
\begin{theorem}\label{thm:ExC3}
  Let $\psi: N^2 \to \mathbb R^3:(v,w) \mapsto (v,w,f(v,w))$ be a
  positive definite improper affine sphere with affine normal
  $(0,0,1)$. Define $\hs(t,v,w)= (t v, t w,t f(v,w) -c t^4,t)$
  or $\hs(t,v,w)= (v, w, f(v,w) +c t^3,t^4)$, where $c,t\in \R^+$.
  Then $\hs$ defines a $3$-dimensional indefinite improper affine
  hypersphere, which admits a pointwise $\mathbb
  Z_3$- or $SO(2)$-symmetry.
\end{theorem}

\section[Pointwise $\SO(1,1)$-symmetry]{Pointwise $\boldsymbol{\SO(1,1)}$-symmetry}
\label{sec:type8}

Let $\M$ be a hypersphere admitting a $\SO(1,1)$-symmetry. We only state the classif\/ication results. The proofs are done quite similar, using a lightvector-frame instead of an orthonormal one, and will appear elsewhere. We denote a lightvector-frame by $\{E,V,F\}$, where $E$ and $F$ are lightvectors and $V$ is spacelike (cf.~\cite{S07a}).

\begin{lemma}\label{lem:KfT8}
Let $\M$ be an affine hypersphere in $\mathbb R^4$ which admits a
  pointwise $\SO(1,1)$-symmetry. Let $p \in M$. Then there exists a
  lightvector-frame $\{E,V,F\}$ defined in a neighborhood of the point~$p$ and a positive function $b_4$ such that $K$ is given by:
\begin{alignat*}{4}
& K(V,V)= -2b_4 V,\qquad & & K(V,E)=b_4 E, \qquad && K(V,F)=
b_4 F,&\\
& K(E,E)=0,\qquad && K(E,F)= b_4 V, \qquad &&
K(F,F)=0.&
\end{alignat*}
\end{lemma}

In the following we denote the coef\/f\/icients of the Levi-Civita
connection with respect to this frame by:
\begin{alignat*}{4}
  & \widehat{\nabla}_E E = \ac11 E + \bc11 V,\qquad &&
   \widehat{\nabla}_E V = \ac12 E - \bc11 F,\qquad &&
   \widehat{\nabla}_E F =-\ac12 V - \ac11 F, & \\
  & \widehat{\nabla}_V E = \ac21 E + \bc21 V, \qquad &&
   \widehat{\nabla}_V V = \ac22 E - \bc21 F, \qquad &&
   \widehat{\nabla}_V F =-\ac22 V - \ac21 F,& \\
  & \widehat{\nabla}_F E = \ac31 E + \bc31 V,\qquad &&
   \widehat{\nabla}_F V = \ac32 E - \bc31 F, \qquad &&
   \widehat{\nabla}_F F =-\ac32 V - \ac31 F.&
 \end{alignat*}

Similar as before, it turns out that the vector f\/ield $V$ is globally def\/ined, and we can def\/ine the
distributions $L_1=\Span\{V\}$ and $L_2=\Span\{E,F\}$. Again, $L_1$ is autoparallel with
respect to~$\widehat\nabla$, and $L_2$ is spherical with mean curvature normal
  $-\ac12 V$. We introduce a coordinate function~$v$ by $\ptv:=V$.

\begin{lemma}\label{lem:nut8}
  The function $\nu=b_4^2-\ac12^2-H$ is globally defined,
  $\ptv(e^{2f} \nu)=0$ and $\nu$ vanishes identically or nowhere on $\R$.
\end{lemma}

Again we have to distinguish three cases.

\subsection[The first case: $\nu \neq 0$ on $\M$]{The f\/irst case: $\boldsymbol{\nu \neq 0}$ on $\boldsymbol{\M}$}
\label{sec:case1t8}

\begin{theorem}\label{thm:ClassC1t8} Let $\M$ be an indefinite affine hypersphere
of $\Rf$ which admits a pointwise $SO(1,1)$-symmetry. Let
$b_4^2-\ac12^2 \neq H$ for some $p\in \M$. Then $\M$ is affine
equivalent to
\[
\hs: \ I\times N^2\to \Rf: \ (v,x,y)\mapsto (\gamma_1(v),
  \gamma_2(v) \phi(x,y)),
  \]
  where
  $\phi: N^2 \to \mathbb R^3$ is a one-sheeted hyperboloid
  and $\gamma:I\to \mathbb R^2$ is a~curve.
\end{theorem}

\begin{theorem}\label{thm:ConstrC1t8}
  Let $\phi:N^2 \to \mathbb R^3$ be a one-sheeted hyperboloid and
  let $\gamma: I \to \mathbb R^2$ be a curve. Define $\hs:I\times N^2\to \Rf$ by
  $\hs(v,x,y)=(\gamma_1(v), \gamma_2(v) \phi(x,y ))$.
  \begin{enumerate}\itemsep=0pt
  \item[$(i)$] If $\ga=(\ga_1,\ga_2)$ satisfies $\ga_2^2|\ga_1\ga_2' - \ga_1'
  \ga_2|^5= |\ga_1'\ga_2'' - \ga_1'' \ga_2'|
  (\ga_1')^2\neq 0$, then $\hs$ defines a~$3$-di\-men\-sional indefinite
  proper affine hypersphere.
  \item[$(ii)$] If $\ga=(\ga_1,\ga_2)$ satisfies $\ga_2^2|\ga_2'|^5=
  |\ga_1'\ga_2'' - \ga_1'' \ga_2'|
  (\ga_1')^2\neq 0$, then $\hs$ defines a $3$-dimensional indefinite
  improper affine hypersphere.
  \end{enumerate}
\end{theorem}

\begin{theorem}\label{thm:ExC1t8}
  Let $\phi:N^2 \to \mathbb R^3$ be a one-sheeted hyperboloid and
  let $\gamma: I \to \mathbb R^2$ be a curve, such that
  $\hs(v,x,y)=(\gamma_1(v), \gamma_2(v) \phi(x,y ))$ defines
  a $3$-dimensional indefinite affine hypersphere. Then
  $\hs(I\times N^2)$ admits a pointwise $SO(1,1)$-symmetry.
\end{theorem}

\subsection[The second case: $\nu \equiv 0$ and $H\neq 0$ on $\M$]{The second case: $\boldsymbol{\nu \equiv 0}$ and $\boldsymbol{H\neq 0}$ on $\boldsymbol{\M}$}
\label{sec:case2t8}

\begin{theorem} \label{thm:ClassC2t8} Let $\M$ be an indefinite proper affine
  hypersphere of $\Rf$ which admits a pointwise
 $SO(1,1)$-symmetry. Let
  $H= b_4^2-\ac12^2 (\neq 0)$ on $\M$. Then $\M$ is affine equivalent with
\[
\hs: \ I\times N^2\to \Rf: \ (v,x,y)\mapsto (\gamma_1(v) x,\gamma_1(v)
  y,\gamma_1(v) f(x,y) + \gamma_2(v), \gamma_1(v)),
  \] where $\psi: N^2
  \to \mathbb R^3:(x,y) \mapsto (x,y,f(x,y))$ is a hyperbolic
  paraboloid with affine normal $(0,0,1)$ and $\gamma:I\to \mathbb
  R^2$ is a curve.
\end{theorem}

\begin{theorem}\label{thm:ConstrC2t8}
  Let $\psi:N^2 \to \mathbb R^3$ be a hyperbolic paraboloid with
  affine normal $(0,0,1)$, and let $\gamma: I \to \mathbb R^2$ be a
  curve. Define $\hs:I\times N^2\to \Rf$ by $\hs(v,x,y)= (\gamma_1(v) x,\gamma_1(v)
  y,\gamma_1(v) f(x,y) + \gamma_2(v), \gamma_1(v))$. If $\ga=(\ga_1,\ga_2)$ satisfies $\ga_1^2|\ga_1\ga_2' -
  \ga_1' \ga_2|^5= |\ga_1'\ga_2'' - \ga_1'' \ga_2'| (\ga_1')^2\neq 0$,
  then $\hs$ defines a $3$-dimensional indefinite proper affine
  hypersphere.
\end{theorem}

\begin{theorem}\label{thm:ExC2t8}
  Let $\psi:N^2 \to \mathbb R^3$ be a hyperbolic paraboloid with
  affine normal $(0,0,1)$, and let $\gamma: I \to \mathbb R^2$ be a
  curve, such that $\hs(v,x,y)= (\gamma_1(v) x,\gamma_1(v)
  y,\gamma_1(v) f(x,y) + \gamma_2(v), \gamma_1(v))$ defines a~$3$-di\-men\-sional indefinite proper affine hypersphere. Then
  $\hs(I\times N^2)$ admits a pointwise $SO(1,1)$-symmetry.
\end{theorem}

\subsection[The third case: $\nu \equiv 0$ and $H=0$ on $\M$]{The third case: $\boldsymbol{\nu \equiv 0}$ and $\boldsymbol{H=0}$ on $\boldsymbol{\M}$}
\label{sec:case3t8}

\begin{theorem} \label{thm:ClassC3t8} Let $\M$ be an indefinite improper affine
  hypersphere of $\Rf$ which admits a pointwise $SO(1,1)$-symmetry.
  Let $\ac12^2 =b_4^2$ on $M^3$. Then $\M$ is affine
  equivalent with either
  \begin{gather*}
  \hs: \ I\times N^2\to \Rf: \ (v,x,y)\mapsto (vx, vy, v f(x,y)
  - c v^4, v),\quad(\ac12 =b_4)\qquad \text{or}\\
   \hs: \ I\times N^2\to \Rf: \ (v,x,y)\mapsto (x,y, f(x,y) + c v^3, v^4 )
  ,\qquad(-\ac12 =b_4),
\end{gather*}
 where $\psi: N^2 \to \mathbb R^3:(v,w) \mapsto (v,w,f(v,w))$ is a
  hyperbolic paraboloid with affine normal
  $(0,0,1)$ and $c,t\in \R^+$.
\end{theorem}

\begin{theorem}\label{thm:ExC3t8}
  Let $\psi: N^2 \to \mathbb R^3:(v,w) \mapsto (v,w,f(v,w))$ be a
  hyperbolic paraboloid with affine nor\-mal
  $(0,0,1)$. Define $\hs(t,v,w)= (t v, t w,t f(v,w) -c t^4,t)$
  or $\hs(t,v,w)= (v, w, f(v,w) +c t^3,t^4)$, where $t\in \R^+$, $c\neq 0$.
  Then $\hs$ defines a $3$-dimensional indefinite improper affine
  hypersphere, which admits a pointwise $SO(1,1)$-symmetry.
\end{theorem}

\appendix

\section[Computations for pointwise $\SO(2)$-, $S_3$- or $\Z_3$-symmetry]{Computations for pointwise $\boldsymbol{\SO(2)}$-, $\boldsymbol{S_3}$- or $\boldsymbol{\Z_3}$-symmetry}

\noindent\(e[1]\text{:=}\{1,0,0\};e[2]\text{:=}\{0,1,0\};e[3]\text{:=}\{0,0,1\}\);

\subsection*{ONB of SO(1,2): \textit{e}[1]=T, \textit{e}[2]=V, \textit{e}[3]=W}

\subsection*{Af\/f\/ine metric}

\noindent h[\text{y$\_$},\text{z$\_$}]\text{:=}-y[[1]]z[[1]]+\text{Sum}[y[[i]]z[[i]],\{i,2,3\}];

\subsection*{Dif\/ference tensor (r=a1, s=a4, u=a6)}

\noindent K[\text{y$\_$},\text{z$\_$}]\text{:=} \text{Sum}[y[[i]]z[[j]]k[i,j], \{i,1,3\}, \{j,1,3\}];\\
 k[1,1]\text{:=}\{-r,0,0\};
 k[1,2]\text{:=}\{0,s,0\};
 k[1,3]\text{:=}\{0,0,r-s\};\\
 k[2,1]\text{:=}\{0,s,0\};
 k[2,2]\text{:=}\{-s,u,0\};
 k[2,3]\text{:=}\{0,0,-u\};\\
 k[3,1]\text{:=}\{0,0,r-s\};
 k[3,2]\text{:=}\{0,0,-u\};
 k[3,3]\text{:=}\{-(r-s),-u,0\};

\subsection*{Ricci tensor, scalar curvature and Pick invariant}

\subsubsection*{[Kx,Ky]z}

\noindent\(\text{LK}[\text{x$\_$},\text{y$\_$},\text{z$\_$}]\text{:=}K[x,K[y,z]]-K[y,K[x,z]];\)

\noindent\(\text{ListLK}\text{:=}\{\text{LK}[e[1],e[2],e[1]],\text{LK}[e[1],e[2],e[2]],\text{LK}[e[1],e[2],e[3]],\\
\text{LK}[e[1],e[3],e[1]],\text{LK}[e[1],e[3],e[2]],\text{LK}[e[1],e[3],e[3]],\\
\text{LK}[e[2],e[3],e[1]],\text{LK}[e[2],e[3],e[2]],\text{LK}[e[2],e[3],e[3]]\};\\
\text{FullSimplify}[\text{ListLK}]\)

\subsubsection*{Curvature tensor (of the Levi-Civita connection)}

\noindent\(R[\text{x$\_$},\text{y$\_$},\text{z$\_$}]\text{:=}H(h[y,z] x-h[x,z] y)-K[x,K[y,z]]+K[y,K[x,z]];\)

\noindent\(\text{ListR}\text{:=}\{R[e[1],e[2],e[1]],R[e[1],e[2],e[2]],R[e[1],e[2],e[3]],\\
R[e[1],e[3],e[1]],R[e[1],e[3],e[2]],R[e[1],e[3],e[3]],\\
R[e[2],e[3],e[1]],R[e[2],e[3],e[2]],R[e[2],e[3],e[3]]\};\\
\text{Simplify}[\text{ListR}]\)

\subsubsection*{Ricci tensor (of the Levi-Civita connection)}

\noindent\(\text{ric}[\text{x$\_$},\text{y$\_$}]\text{:=}
\text{Simplify}[
(-h[R[e[1],x,y],e[1]]+h[R[e[2],x,y],e[2]]+
h[R[e[3],x,y],e[3]])]\);

\noindent\(\text{Listric}\text{:=}\{\text{ric}[e[1],e[1]],\text{ric}[e[1],e[2]],\text{ric}[e[1],e[3]],
\text{ric}[e[2],e[2]],\text{ric}[e[2],e[3]],\text{ric}[e[3],e[3]]\}; \\
\text{Simplify}[\text{Listric}]\)

\subsubsection*{Scalar curvature (of the Levi-Civita connection)}

\noindent\(\text{sc}\text{:=}1/6(-\text{ric}[e[1],e[1]]+\text{ric}[e[2],e[2]]+\text{ric}[e[3],e[3]]); \\
\text{Simplify}[\text{sc}]\)

\subsubsection*{Pick invariant}

\noindent\(P\text{:=}
1/6(-(k[1,1][[1]]){}^{\wedge}2 + (k[2,2][[2]]){}^{\wedge}2 +(k[3,3][[3]]){}^{\wedge}2 +
3 ((k[1,1][[2]]){}^{\wedge}2+(k[1,1][[3]]){}^{\wedge}2 -(k[2,2][[1]]){}^{\wedge}2 -
(k[3,3][[1]]){}^{\wedge}2 +(k[2,2][[3]]){}^{\wedge}2+(k[3,3][[2]]){}^{\wedge}2)-
6 (k[1,2][[3]]){}^{\wedge}2); \\
\text{Simplify}[P]\)

\subsection*{Lemma 1}

\noindent\(r=2s; \ \text{Simplify}[\text{Listric}]\)

\noindent\(\text{Simplify}[P]\)

\subsection*{Lemma 8}

\noindent\(R[e[2],e[3],e[2]]\)

\subsection*{Lemma 3 }

\subsubsection*{Levi-Civita connection (ONB)}

\noindent\(n[1,1]\text{:=}\{0,\text{a12},\text{a13}\};
n[1,2]\text{:=}\{\text{a12},0,-\text{b13}\};
n[1,3]\text{:=}\{\text{a13},\text{b13},0\};\\
n[2,1]\text{:=}\{0,\text{a22},\text{a23}\};
n[2,2]\text{:=}\{\text{a22},0,-\text{b23}\};
n[2,3]\text{:=}\{\text{a23},\text{b23},0\};\\
n[3,1]\text{:=}\{0,\text{a32},\text{a33}\};
n[3,2]\text{:=}\{\text{a32},0,-\text{b33}\};
n[3,3]\text{:=}\{\text{a33},\text{b33},0\};\)

\noindent\(\text{Na}[\text{y$\_$},\text{z$\_$}]\text{:=}\text{Sum}[y[[i]]z[[j]]n[i,j],\{i,1,3\},\{j,1,3\}]+
\text{Sum}[y[[1]]\text{Dt}[z[[i]],\text{f1}]e[i]\)\\
\(+\, y[[2]]\text{Dt}[z[[i]],\text{f2}]e[i] +
y[[3]]\text{Dt}[z[[i]],\text{f3}]e[i],\{i,1,3\}];\)

\subsubsection*{Codazzi for K (af\/f\/ine hypersphere)}

\noindent
\(\text{codazzi}[\text{x$\_$},\text{y$\_$},\text{z$\_$}]\text{:=}\text{Na}[x,K[y,z]]-K[\text{Na}[x,y],z]-K[y,\text{Na}[x,z]] -
\text{Na}[y,K[x,z]]+K[\text{Na}[y,x],z]+K[x,\text{Na}[y,z]] ;\)

\noindent\(\text{eq1}\text{:=}\text{Simplify}[\text{codazzi}[e[2],e[1],e[1]]];
\text{eq2}\text{:=}\text{Simplify}[\text{codazzi}[e[3],e[1],e[1]]];\\
\text{eq3}\text{:=}\text{Simplify}[\text{codazzi}[e[1],e[2],e[2]]];
\text{eq4}\text{:=}\text{Simplify}[\text{codazzi}[e[3],e[2],e[2]]];\\
\text{eq5}\text{:=}\text{Simplify}[\text{codazzi}[e[1],e[3],e[3]]];
\text{eq6}\text{:=}\text{Simplify}[\text{codazzi}[e[2],e[3],e[3]]];\\
\text{eq7}\text{:=}\text{Simplify}[\text{codazzi}[e[1],e[2],e[3]]];
\text{eq8}\text{:=}\text{Simplify}[\text{codazzi}[e[2],e[3],e[1]]];\\
\text{eq9}\text{:=}\text{Simplify}[\text{codazzi}[e[3],e[1],e[2]]];\)

\noindent\(\text{eq}\text{:=}\{\text{eq1},\text{eq2},\text{eq3},\text{eq4},\text{eq5},\text{eq6},\text{eq7},\text{eq8},\text{eq9}\}; \text{eq}\)

\subsubsection*{1. case: u${}^{\wedge}$2$\neq $4 s${}^{\wedge}$2}

{\bf conclusions from eq1,2,4:}

\noindent\(\text{Simplify}[\text{eq2}[[1]]-2 \text{eq4}[[1]]]\)

\noindent\(\text{Simplify}[\text{eq1}[[3]]+ \text{eq2}[[2]]]\)

\noindent\(\text{a13}=0; \text{a32}=-\text{a23};\)

\noindent\(\text{eq}\)

\noindent {\bf conclusions from eq1,2,6:}

\noindent\(\text{Simplify}[-2 \text{eq6}[[1]]+ \text{eq1}[[1]]]\)

\noindent\(\text{Simplify}[\text{eq1}[[2]]- \text{eq2}[[3]]]\)

\noindent\(\text{a12}=0; \text{a33}=\text{a22};\)

\noindent\(\text{eq}\)

\noindent\(\text{Clear}[\text{a13},\text{a12},\text{a32},\text{a33}]\)

\subsubsection*{2. case: u=2s$\neq $0}

\noindent\(u=2s;\ \text{eq}\)

\noindent {\bf conclusions from eq8, eq1:}

\noindent\(\text{a32}= \text{a23}; \text{a13}=- 2 \text{a23};\text{eq}\)

\noindent {\bf conclusions from eq3:}

\noindent\(\text{b13}=0; \text{eq}\)

\noindent\(\text{Simplify}[ \text{eq1}[[2]]- \text{eq2}[[3]]]\)

\noindent\(\text{a12}=-(\text{a33}- \text{a22}); \text{eq}\)

\noindent\(\text{Simplify}[ \text{eq3}[[2]]-1/2 \text{eq1}[[1]]+2 \text{eq1}[[2]]]\)

\noindent\(\text{a33}=-\text{a22}; \text{eq}\)

 It follows that T(a4)=0, V(a4)=-4a22 a4, W(a4)=4a23 a4.

\noindent\(\text{Simplify}[ \text{eq4}[[2]]+\text{eq2}[[1]]]\)

\noindent\(\text{Simplify}[ \text{eq4}[[3]]+ \text{eq1}[[1]]]\)

\noindent\(\text{b23}=- \text{a23}; \text{b33}=- \text{a22};\)

\subsection*{Lemma 4}
\subsubsection*{Gauss for Levi-Civita connection (af\/f\/ine hypersphere)}

\noindent\(\text{gaussLC}[\text{x$\_$},\text{y$\_$},\text{z$\_$}]\text{:=}\text{Na}[x,\text{Na}[y,z]]-\text{Na}[y,\text{Na}[x,z]]-
\text{Na}[\text{Na}[x,y]-\text{Na}[y,x],z]-H h[y,z]x+H h[x,z]y+K[x,K[y,z]]-
K[y,K[x,z]];\)

\noindent\(\text{eq11}\text{:=}\text{Simplify}[\text{gaussLC}[e[1],e[2],e[2]]];
\text{eq12}\text{:=}\text{Simplify}[\text{gaussLC}[e[1],e[3],e[2]]];\\
\text{eq13}\text{:=}\text{Simplify}[\text{gaussLC}[e[2],e[3],e[2]]];
\text{eq14}\text{:=}\text{Simplify}[\text{gaussLC}[e[1],e[2],e[1]]];\\
\text{eq15}\text{:=}\text{Simplify}[\text{gaussLC}[e[1],e[3],e[1]]];
\text{eq16}\text{:=}\text{Simplify}[\text{gaussLC}[e[2],e[3],e[1]]];\\
\text{eq17}\text{:=}\text{Simplify}[\text{gaussLC}[e[1],e[2],e[3]]];
\text{eq18}\text{:=}\text{Simplify}[\text{gaussLC}[e[1],e[3],e[3]]];\\
\text{eq19}\text{:=}\text{Simplify}[\text{gaussLC}[e[2],e[3],e[3]]];\)

\noindent\(\text{eqG}\text{:=}\{\text{eq11},\text{eq12},\text{eq13},\text{eq14},\text{eq15},\text{eq16},\text{eq17},\text{eq18},\text{eq19}\};
\)

\subsubsection*{2. case: u=2s $\neq $0}

\noindent\(\text{eqG}\)

\noindent\(\text{Simplify}[\text{eq11}[[1]]-\text{eq12}[[3]]]\)

\noindent\(\text{Simplify}[\text{eq15}[[3]]+\text{eq12}[[3]]]\)

Contradiction to eq13.3

\noindent\(\text{Clear}[\text{b33},\text{b23},\text{a33},\text{a12},\text{b13},\text{a32},\text{a13},u]\)

\subsubsection*{1. case: u${}^{\wedge}$2$\neq $4 s${}^{\wedge}$2}

\noindent\(\text{a13}=0; \text{a32}=-\text{a23};\text{a12}=0; \text{a33}=\text{a22};\text{eqG}\)

\subsection*{Acknowledgements}

Partially supported by the DFG-Project PI 158/4-5
`Geometric Problems and Special PDEs'.

\pdfbookmark[1]{References}{ref}
\LastPageEnding

\end{document}